\documentclass[12pt]{article}%
\usepackage{amssymb}
\usepackage[author-year]{amsrefs}
\usepackage{amsmath}
\usepackage{amsfonts}
\usepackage{float}
\usepackage{graphicx}
\usepackage{amscd}
\usepackage{amsthm}
\usepackage{amsxtra}
\usepackage{verbatim}
\usepackage{endnotes}
\usepackage{times}
\usepackage{hyperref}%
\theoremstyle{plain}
\newtheorem{X}{X}[section]
\newtheorem{theorem}[X]{Theorem}
\newtheorem{proposition}[X]{Proposition}
\newtheorem{lemma}[X]{Lemma}
\newtheorem{corollary}[X]{Corollary}

\theoremstyle{definition}
\newtheorem{definition}[X]{Definition}
\newtheorem{example}[X]{Example}
\newtheorem{plain}[X]{}

\newtheorem{remark}[X]{Remark}
\newtheorem{exercise}[X]{Exercise}

\newcounter{junk}

\newenvironment{mylist}
{\begin{list}
{--}
{\setlength{\leftmargin}{.5in}\setlength{\rightmargin}{.5in}}}
{\end{list}}

\renewenvironment{itemize}
{\begin{mylist}}
{\end{mylist}}

\def\1{{1\mkern-7mu1}}  

\newcommand\Aut{\operatorname{Aut}}

\newcommand\Br{\operatorname{Br}}

\newcommand\dlim{\varinjlim}

\newcommand\Gal{\operatorname{Gal}}

\newcommand\Hom{\operatorname{Hom}}

\newcommand\id{\operatorname{id}}

\newcommand\im{\operatorname{Im}}  

\newcommand\inv{\operatorname{inv}}

\newcommand\Ker{\operatorname{Ker}}

\newcommand\Nm{\operatorname{Nm}}

\newcommand\ob{\operatorname{ob}}
\newcommand\ord{\operatorname{ord}}

\newcommand\plim{\varprojlim}

\newcommand\Res{\operatorname{Res}}

\newcommand\Sh{\operatorname{Sh}}

\newcommand\Spec{\operatorname{Spec}}

\def\Aff{{\mathsf{Aff}}}
\def\CM{{\mathsf{CM}}}

\def\Mot{{\mathsf{Mot}}}

\def\Rep{{\mathsf{Rep}}}

\def\Vc{{\mathsf{Vec}}}

\def\ord{{\mathrm{ord}}}

\def\Tors{{\textsc{Tors}}}

\let\xr=\xrightarrow

\makeatletter

\let\citeright\@empty
\def\citesel@object#1#2#3#4#5{\PrintCiteNames{#3} \@citeleft#4}
\makeatother
\renewcommand{\cite}{\ocite}

\begin{document}

\title{Gerbes and Abelian Motives}
\date{\today. v1.}
\author{J.S. Milne\thanks{Partially supported by the National Science Foundation.}}
\maketitle

\begin{abstract}
Assuming the Hodge conjecture for abelian varieties of CM-type, one obtains a
good category of abelian motives over $\mathbb{F}_{p}^{\text{al}}{}$ and a
reduction functor to it from the category of CM-motives over $\mathbb{Q}%
{}^{\text{al}}$. Consequently, one obtains a morphism${}$ of gerbes of fibre
functors with certain properties. We prove unconditionally that there exists a
morphism of gerbes with these properties, and we classify them.

\end{abstract}
\tableofcontents

\section*{Introduction}

Fix a $p$-adic prime of $\mathbb{Q}^{\text{al}}$, and denote its residue field
by $\mathbb{F}{}$. Let $\CM(\mathbb{Q}^{\text{al}})$ be the category of
motives based on the abelian varieties of CM-type over $\mathbb{Q}^{\text{al}%
}$, and let $\Mot(\mathbb{F}{})$ be the similar category based on the abelian
varieties over $\mathbb{F}{}$; in both cases, the correspondences are to be
the algebraic cycles modulo numerical equivalence. Both are Tannakian
$\mathbb{Q}{}$-categories and, because every abelian variety of CM-type has
good reduction, there is an exact tensor \textquotedblleft
reduction\textquotedblright\ functor%
\[
R\colon\CM(\mathbb{Q}^{\text{al}})\rightarrow\Mot(\mathbb{F}{})\text{.}%
\]
From $R$, we obtain a morphism%
\[
R^{\vee}\colon\Mot(\mathbb{F}{})^{\vee}\rightarrow\CM(\mathbb{Q}^{\text{al}%
})^{\vee}%
\]
of gerbes of fibre functors. The bands of ${}$ $\Mot(\mathbb{F}{})^{\vee}$ and
$\CM(\mathbb{Q}^{\text{al}})^{\vee}{}$ are commutative, and so $R^{\vee}{}$ is
bound by a homomorphism
\[
\rho\colon P\rightarrow S
\]
of commutative affine group schemes over $\mathbb{Q}{}$.

It is known how to construct an explicit homomorphism $\rho\colon P\rightarrow
S$ of commutative affine group schemes (Grothendieck, Langlands,
Rapoport,\ldots) and to prove that it becomes the homomorphism in the last
paragraph when the Hodge conjecture is assumed for abelian varieties of
CM-type (this is explained in \cite{milne1994a}, \S 2, \S 4, under the
additional assumption of the Tate conjecture for abelian varieties over finite
fields, but that assumption is shown to be superfluous in \cite{milne1999b}).

In this paper, we construct (unconditionally) an explicit morphism\qquad%
\[
\mathsf{P{}}\rightarrow\mathsf{S{}}%
\]
of gerbes bound by the homomorphism $\rho\colon P\rightarrow S$ in the last
paragraph, and we prove that it becomes the motivic morphism in the first
paragraph when the Hodge conjecture is assumed for abelian varieties of
CM-type. Moreover, we classify the morphisms $\mathsf{P{}}\rightarrow
\mathsf{S{}}$ bound by $\rho$ and satisfying certain natural conditions.

Upon choosing an $x\in\mathsf{P}{}(\mathbb{Q}{}^{\text{al}})$, we obtain a
morphism of groupoids%
\[
\mathsf{P{}}(x)\rightarrow\mathsf{S{}}(x)
\]
such that the map on kernels%
\[
\mathsf{P{}}(x)^{\Delta}\rightarrow\mathsf{S{}}(x)^{\Delta}%
\]
is the homomorphism $\rho$.

We now describe the contents of the paper in more detail.

In \S 1 we compute some inverse limits. In particular, we show that certain
inverse limits are infinite (see, for example, \ref{il06}), and hence can not
be ignored as they have been in previous works.

After some preliminaries on the cohomology of protori in \S 2, we construct
and classify certain \textquotedblleft fundamental\textquotedblright%
\ cohomology classes in \S \S 3--4.

In \S 5, we review part of the theory of gerbes and their classification
(\cite{giraud1971}, \cite{debremaeker1977a}).

In \S 6, we prove the existence of gerbes $\mathsf{P}$ having the properties
expected of $\Mot(\mathbb{F}{})^{\vee}$, and we classify them.

Finally, in \S 7 we prove the existence of morphism of gerbes $\mathsf{P}%
\rightarrow\mathsf{S}$ having the properties expected of $\Mot(\mathbb{F}%
{})^{\vee}\rightarrow\CM(\mathbb{Q}^{\text{al}})^{\vee}$, and we classify them.

In large part, this article is a critical re-examination of the results and
proofs in \S \S 2--4 (pp118--165) of \cite{langlands1987}. In particular, we
eliminate the confusion between fpqc cohomology groups and inverse limits of
Galois cohomology groups that has persisted in the literature for fifteen
years (e.g., \cite{langlands1987}, \cite{milne1994a}, \cite{reimann1997}),
which amounted to setting certain $\varprojlim^{1}$s equal to zero. Also we
avoid the confusion between gerbes and groupoids to be found in
\cite{langlands1987}. Finally, we are concerned, not just with the existence
of the various object, but also to what extent they are well defined. See the
notes at the end of the sections for more details.

For the convenience of the reader, I have made this article independent of
earlier articles on the topic.

\subsection{Notations and conventions}

The algebraic closure of $\mathbb{Q}{}$ in $\mathbb{C}{}$ is denoted
$\mathbb{Q}{}^{\text{al}}$. The symbol $p$ denotes a fixed finite prime of
$\mathbb{Q}{}$ and $\infty$ denotes the real prime.

Complex conjugation on $\mathbb{C}{}$, or a subfield, is denoted by $\iota$. A
CM-field is a finite extension $E$ of $\mathbb{Q}{}$ admitting an involution
$\iota_{E}\neq1$ such that $\rho\circ\iota_{E}=\iota\circ\rho$ for all
homomorphisms $\rho\colon E\rightarrow\mathbb{C}{}$. The composite of the
CM-subfields of $\mathbb{Q}{}^{\text{al}}$ is denoted $\mathbb{Q}{}%
^{\text{cm}}$.

The set $\mathbb{N}{}\smallsetminus\{0\}$, partially ordered by divisibility,
is denoted $\mathbb{N}{}^{\times}$. For a finite set $S$, $\mathbb{Z}{}%
^{S}=\Hom(S,\mathbb{Z}{})$.

For a perfect field $k$ of characteristic $p$, $W(k)$ is the ring of Witt
vectors with coefficients in $k$ and $B(k)$ is the field of fractions of
$W(k)$. The automorphism of $W(k)$ inducing $x\mapsto x^{p}$ on the residue
field is denoted $\sigma$.

For a finite extension of fields $K\supset k$, $(\mathbb{G}_{m})_{K/k}$ is the
torus over $k$ obtained from $\mathbb{G}_{m}$ over $K$ by restriction of scalars.

Groupoid will always mean affine transitive groupoid scheme
(\cite{deligne1990}).

For a group scheme $G$, a right $G$-object $X$, and a left $G$-object $Y$,
$X\wedge^{G}Y$ denotes the contracted product of $X$ and $Y$, that is, the
quotient of $X\times Y$ by the diagonal action of $G$, $(x,y)g=(xg,g^{-1}y)$.
When $G\rightarrow H$ is a homomorphism of group schemes, $X\wedge^{G}H$ is
the $H$-object obtained from $X$ by extension of the structure group. In this
last case, if $X$ is a $G$-torsor, then $X\wedge^{G}H$ is also an $H$-torsor.

The notation $X\approx Y$ means that $X$ and $Y$ are isomorphic, and $X\cong
Y$ means that $X$ and $Y$ are canonically isomorphic (or that a particular
isomorphism is given). 

Direct and inverse limits are always with respect to directed index sets
(partially ordered set $I$ such that, for all $i,j\in I$, there exists a $k\in
K$ for which $i\leq k$ and $j\leq k$). An inverse system is \emph{strict }if
its transition maps are surjective.

For class field theory we use the sign convention that the local Artin map
sends a prime element to the Frobenius automorphism (that inducing $x\mapsto
x^{q}$ on the residue field); for Hodge structures we use the convention that
$h(z)$ acts on $V^{p,q}$ as $z^{p}\bar{z}^{q}$.

\subsection{Philosophy\label{philosophy}}

We adopt the following definitions.

\begin{plain}
\label{p1}An element of a set is \emph{well-defined} (by a property,
construction, condition, etc.) when it is uniquely determined (by the
property, construction, condition etc.).
\end{plain}

\begin{plain}
\label{p2}An object of a category is \emph{well-defined} when it is uniquely
determined up to a uniquely given isomorphism. For example, an object defined
by a universal property is well-defined in this sense. When $X$ is
well-defined in this sense and $X^{\prime}$ also has the defining property,
then each element of $X$ corresponds to a well-defined element of $X^{\prime}$.
\end{plain}

\begin{plain}
\label{p3}A category is \emph{well-defined} when it is uniquely determined up
to a category-equivalence\footnote{Recall that a category equivalence is a
functor $F\colon\mathsf{C}\rightarrow\mathsf{C}^{\prime}$ for which there
exists a functor $G\colon\mathsf{C}^{\prime}\rightarrow\mathsf{C}$ and
isomorphisms $\phi\colon\id\rightarrow GF$ and $\psi\colon1\rightarrow FG$
such that $F\phi=\psi F$ (\cite{bucur1968}, 1.6; $F\phi$ and $\psi F$ are
morphisms $F\rightarrow FGF$). Our condition says that if $\mathsf{C}$ and
$\mathsf{C}^{\prime}$ have the defining property, then there is a
distinguished class of equivalences $F\colon\mathsf{C}\rightarrow
\mathsf{C}^{\prime}$ and a distinguished class of equivalences $F^{\prime
}\colon\mathsf{C}^{\prime}\rightarrow\mathsf{C}$; if $F_{1}\colon
\mathsf{C}\rightarrow\mathsf{C}^{\prime}$ and $F_{2}\colon\mathsf{C}%
\rightarrow\mathsf{C}^{\prime}$ are both in the distinguished class, then
there is given a (unique) isomorphism $F_{1}\rightarrow F_{2}$.} that is
itself uniquely determined up to a uniquely given isomorphism. When
$\mathsf{C}{}$ is well-defined in this sense and $\mathsf{C}{}^{\prime}$ also
has the defining property, then each object of $\mathsf{C}{}$ corresponds to a
well-defined object of $\mathsf{C}{}^{\prime}$.
\end{plain}

\subsection{Advice to the reader}

The article has been written in logical order. The reader is advised to begin
with \S \S 5,6,7 and refer back to the earlier sections only as needed. Also,
the results in \S 6 on the gerbe conjecturally attached to $\Mot(\mathbb{F}%
{})$ require very little of the rest of the paper.

\newpage

\section{Some inverse limits.}

We compute some $\varprojlim$s and $\varprojlim^{1}$s that are needed in the
rest of the paper.

\subsection{Review of higher inverse limits}

{}For an inverse system $(A_{n},u_{n})$ of abelian groups indexed by
$(\mathbb{N},\leq)$, $\varprojlim A_{n}$ and $\varprojlim^{1}A_{n}$ can be
defined to be the kernel and cokernel respectively of%
\begin{equation}
(\ldots,a_{n},\ldots)\mapsto(\ldots,a_{n}-u_{n+1}(a_{n+1}),\ldots)\colon%
{\textstyle\prod\nolimits_{n}}
A_{n}\xr{1-u}%
{\textstyle\prod\nolimits_{n}}
A_{n}\text{.} \label{e10}%
\end{equation}
From the snake lemma, we see that an inverse system of short exact sequences%
\[
0\rightarrow(A_{n})\rightarrow(B_{n})\rightarrow(C_{n})\rightarrow0
\]
gives rise to an exact sequence%
\begin{equation}
0\rightarrow\varprojlim A_{n}\rightarrow\varprojlim B_{n}\rightarrow
\varprojlim C_{n}\rightarrow\varprojlim{}^{1}A_{n}\rightarrow\varprojlim{}%
^{1}B_{n}\rightarrow\varprojlim{}^{1}C_{n}\rightarrow0. \label{e1}%
\end{equation}
In particular, $\varprojlim$ is left exact and $\varprojlim^{1}$ is right exact.

Recall that an inverse system $(A_{n},u_{n})_{n\in\mathbb{N}{}}$ is said to
satisfy the condition (ML) if, for each $n$, the decreasing chain in $A_{n}$
of the images of the $A_{i}$ for $i\geq n$ is eventually constant.

\begin{proposition}
\label{il01b}The group $\plim{}_{n\in\mathbb{N}{}}^{1}A_{n}=0$ if

\begin{enumerate}
\item the $A_{n}$ are compact and the transition maps are continuous, or

\item the system $(A_{n})_{n\in\mathbb{N}{}}$ satisfies (ML).
\end{enumerate}
\end{proposition}

\begin{proof}
Standard.\footnote{Consider the map%
\[
(a_{0},\ldots,a_{n},\ldots,a_{N+1})\mapsto(\ldots,a_{n}-u_{n+1}(a_{n+1}%
),\ldots)\colon%
{\textstyle\prod\nolimits_{0\leq n\leq N+1}}
A_{n}\rightarrow%
{\textstyle\prod\nolimits_{0\leq n\leq N}}
A_{n}.
\]
Let $x=(x_{n})_{n\in\mathbb{N}{}}\in\prod_{n\in\mathbb{N}{}}A_{n}$, and let
$P_{N}$ be the inverse image of $(x_{n})_{0\leq n\leq N}$ in $%
{\textstyle\prod\nolimits_{0\leq n\leq N+1}}
A_{n}$. We have to show that $\plim{}P_{N}$ is nonempty. The projection
$(a_{n})_{0\leq n\leq N+1}\mapsto a_{N+1}\colon P_{N}\rightarrow A_{N+1}$ is
bijective. In case (a), the $P_{N}$ are compact, and so this case follows from
\cite{bourbaki1989}, I \S 9.6, Proposition 8. In case (b), $(P_{N}%
)_{N\in\mathbb{N}{}}$ satisfies (ML). Let $Q_{N}=\cap_{i}\im(P_{N+i}%
\rightarrow P_{N})$. Then each $Q_{N}$ is nonempty, and $(Q_{N})_{N\in
\mathbb{N}{}}$ is a strict inverse system. Hence $\varprojlim Q_{N}$ is
(obviously) nonempty. As $\varprojlim Q_{N}=\varprojlim P_{N}$, this proves
(b).}
\end{proof}

\begin{remark}
\label{il01e}Consider an inverse system $(A_{n})_{n\in\mathbb{N}{}}$ of finite
groups. If $\varprojlim A_{n}$ is infinite, then it is uncountable. In proving
this we may assume\footnote{As in the preceding footnote.} that the transition
maps are surjective. Because $\varprojlim A_{n}$ is infinite, the orders of
the $A_{n}$ are unbounded, and the Cantor diagonalization process can be applied.
\end{remark}

\begin{remark}
\label{il01g}We shall frequently make use of the following (obvious)
criterion:%
\[
\varprojlim A_{n}=0\text{ if, for all }n\text{, }\cap_{i}\im(A_{n+i}%
\rightarrow A_{n})=0.
\]

\end{remark}

\begin{proposition}
\label{il01f}If the inverse system $(A_{n})_{n\in\mathbb{N}{}}$ fails (ML) and
the $A_{n}$ are countable, then $\varprojlim^{1}A_{n}$ is uncountable.
\end{proposition}

\begin{proof}
See, for example, \cite{milne2002c}, 1.1b.
\end{proof}

\begin{example}
\label{il01m}Let $(G_{n})_{n\in\mathbb{N}{}}$ be an inverse system of
commutative affine algebraic groups over $\mathbb{Q}{}$ with surjective
transition maps. The inverse system $(G_{n}(\mathbb{Q}{}))_{n\in\mathbb{N}{}}$
will usually fail (ML), and so $\varprojlim^{1}G_{n}(\mathbb{Q}{})$ will
usually be uncountable. If $G_{n}$ contains no $\mathbb{Q}{}$-split torus,
then $G_{n}(\mathbb{A}{})/G_{n}(\mathbb{Q}{})$ is compact (\cite{platonov1994}%
, Theorem 5.5, p260), and so $\plim{}^{1}G_{n}(\mathbb{A}{})/G_{n}%
(\mathbb{Q}{})=0$; thus
\[
\plim{}^{1}G_{n}(\mathbb{Q})\rightarrow\plim{}^{1}G_{n}(\mathbb{A}{})
\]
is surjective.
\end{example}

\begin{remark}
\label{il01n}For an arbitrary directed set $I$, the category of inverse
systems of abelian groups indexed by $I$ has enough injectives, and so there
are right derived functors $\varprojlim_{I}^{i}$ of $\varprojlim_{I}$
(\cite{jensen1972} , \S 1). For $I=(\mathbb{N}{},\leq)$, they agree with those
defined above (ibid. p13). If $J$ is a cofinal subset of $J$, then
$\varprojlim_{J}^{i}A_{\alpha}\cong\varprojlim_{I}^{i}A_{\alpha}$. This is
obvious for $i=0$ and it follows for a general $i$ by the usual
derived-functor argument (cf. ibid. 1.9 and statement p12). If $(I,\leq)$ is a
countable directed set, then $I$ contains a cofinal subset isomorphic to
$(\mathbb{N}{},\leq)$ or to a finite segment of $(\mathbb{N}{},\leq)$. In the
first case, the above statements apply to inverse systems indexed by $I$.
\end{remark}

\begin{exercise}
\label{il01p}For an abelian group $A$, let $(A,m)$ denote the inverse system
indexed by $\mathbb{N}{}^{\times}$%
\[
\cdots\leftarrow\underset{n}{A}\,\overset{m}{\leftarrow}\,\underset{mn}%
{A}\,\leftarrow\cdots.
\]
Show:

\begin{enumerate}
\item $\varprojlim(A,m)=0=\varprojlim^{1}(A,m)$ if $NA=0$ for some integer $N$;

\item $\varprojlim(\mathbb{Z}{},m)=0$, $\varprojlim^{1}(\mathbb{Z}{}%
,m)\cong\mathbb{\hat{Z}/\mathbb{Z}{}}$;

\item $\varprojlim(\mathbb{Q}{}/\mathbb{Z}{},m)\cong\mathbb{A}{}_{f}$,
$\varprojlim^{1}(\mathbb{Q}/\mathbb{Z}{},m)=0$.$\quad$[Hint: First show that
$\varprojlim(\mathbb{Q}{}_{l}/\mathbb{Z}{}_{l},m)\cong\mathbb{Q}_{l}$.]
\end{enumerate}
\end{exercise}

\subsection{Notations from class field theory}

For a number field $L$, $U(L)$ is the group of units in the ring of integers
of $L$, $I(L)$ is the group of fractional ideals, and $\mathbb{I}{}(L)$ is the
group of id\`{e}les. For a modulus $\mathfrak{m}{}$ of $L$, $C_{\mathfrak{m}%
}(L)$ is the ray class group, and $H_{\mathfrak{m}}(L)$ is the corresponding
ray class field. When $\mathfrak{m}=1$, we omit it from the notation. Thus,
$C(L)$ is the usual ideal class group and $H(L)$ is the Hilbert class field.

For a finite extension $L/F$ of number fields and moduli $\mathfrak{m}{}$ of
$F$ and $\mathfrak{n}{}$ of $L$ such that $\mathfrak{m}{}|\Nm_{L/F}%
(\mathfrak{n}{})$, there is a commutative diagram%
\begin{equation}
\begin{CD} C_{\mathfrak{n}}(L)@>{\cong}>> \Gal(H_{\mathfrak{n}}(L)/L) \\ @VV{\Nm_{L/F}}V@VV{\sigma \mapsto \sigma |H_{\mathfrak{m}}(F)}V \\ C_{\mathfrak{m}}(F)@>\cong>>\Gal(H_{\mathfrak{m}}(F)/F)\end{CD} \label{1d}%
\end{equation}
(\cite{tate1967}, 3.2, p166). The horizontal isomorphisms are defined by the
Artin symbols $(-,H_{\mathfrak{n}}(L)/L)$ and $(-,H_{\mathfrak{m}{}}(F)/F)$.

\subsection{Inverse systems indexed by fields}

\begin{plain}
\label{il01}Let $\mathcal{F}{}$ be an infinite set of subfields of
$\mathbb{Q}^{\text{al}}$ such that

\begin{enumerate}
\item each field in $\mathcal{F}{}$ is of finite degree over $\mathbb{Q}{}$;

\item the composite of any two fields in $\mathcal{F}{}$ is again in
$\mathcal{F}$.
\end{enumerate}

\noindent Then $(\mathcal{F}{},\subset)$ is a directed set admitting a cofinal
subset isomorphic to $(\mathbb{N}{},\leq)$. For such a set $\mathcal{F}{}$,
the norm maps define inverse systems $(F^{\times})_{F\in\mathcal{F}{}}$,
$(I(F))_{F\in\mathcal{F}{}}$, etc..
\end{plain}

\begin{example}
\label{il05}(a) The set of all totally real fields in $\mathbb{Q}{}%
^{\text{al}}$ of finite degree over $\mathbb{Q}{}$ satisfies (\ref{il01}a,b).
If $f\in\mathbb{R}{}[X]$ is monic with $\deg(f)$ distinct real roots, then any
monic polynomial in $\mathbb{R}{}[X]$ sufficiently close to $f$ will have the
same property (because any quadratic factor of it will have discriminant
$>0$). From this it is easy to construct many totally real extensions of
$\mathbb{Q}{}$. For example, Krasner's lemma shows that, for any finite
extension $L$ of $\mathbb{Q}_{p}$, there exists a finite totally real
extension $F$ of $\mathbb{Q}$ contained in $L$ such that $[F\colon
\mathbb{Q}]=[L\colon\mathbb{Q}_{p}]$ and $F\cdot\mathbb{Q}_{p}=L$. Also, a
standard argument in Galois theory can be modified to show that, for all $n$,
there exist totally real fields with Galois group $S_{n}$ over $\mathbb{Q}{}$
(\cite{wei1993}, 1.6.7).

(b) The set of all CM fields in $\mathbb{Q}{}^{\text{al}}$ of finite degree
over $\mathbb{Q}{}$ satisfies (\ref{il01}a,b). Let $d\in\mathbb{Z}{}$, $d>0$.
For any totally real field $F$, $F[\sqrt{-d}]$ is a CM-field, and every
CM-field containing $\sqrt{-d}$ is of this form. As $F$ runs over the totally
real subfields of $\mathbb{Q}{}^{\text{al}}$, $F[\sqrt{-d}]$ runs over a
cofinal set of CM-subfields of $\mathbb{Q}{}^{\text{al}}$.
\end{example}

\subsection{The inverse system $(C(K))$}

\begin{proposition}
\label{il05m}Let $\mathcal{F}$ be the set of totally real fields in
$\mathbb{Q}{}^{\text{al}}$ of finite degree over $\mathbb{Q}{}$. Then
$\plim{}_{\mathcal{F}{}}C(F)=0$.
\end{proposition}

\begin{proof}
The Hilbert class field $H(F)$ of a totally real field $F$ is again totally
real, and diagram (\ref{1d}) shows that, for any totally real field $L$
containing $H(F)$, the norm map $C(L)\rightarrow C(F)$ is zero. Therefore,
$\plim{}C(F)=0$ by (\ref{il01g}).
\end{proof}

\begin{proposition}
\label{il06}Let $\mathcal{F}$ be the set of CM-fields in $\mathbb{Q}%
{}^{\text{al}}$ of finite degree over $\mathbb{Q}{}$. Then
$\plim{}_{\mathcal{F}{}}C(K)$ is uncountable.
\end{proposition}

Let $K$ be a CM-field with largest real subfield $K_{+}$. Let $\mathfrak{m}{}$
be a modulus for $K$ such that $\iota_{K}\mathfrak{m}{}=\mathfrak{m}{}$ and
the exponent of any finite prime ramified in $K/K_{+}$ is even. Then
$\mathfrak{m}{}$ is the extension to $K$ of a unique modulus $\mathfrak{m}{}%
{}_{+}$ for $K_{+}$ such that $(\mathfrak{m}{}_{+})_{\infty}=1$. This last
condition implies that $H_{\mathfrak{m}{}_{+}}(K_{+})\cap K=K_{+}$, and so the
norm map $C_{\mathfrak{m}{}}(K)\rightarrow C_{\mathfrak{m}{}_{+}}(K_{+})$ is
surjective (apply (\ref{1d})). Define $C_{\mathfrak{m}{}}^{-}(K)$ by the exact
sequence%
\begin{equation}
\begin{CD} 0\rightarrow C_{\mathfrak{m}}^{-}(K)@>>> C_{\mathfrak{m}}(K)@>{\Nm_{K/K_+}}>>C_{\mathfrak{m}_{+}}(K_{+})\rightarrow0\text{.} \end{CD} \label{e3}%
\end{equation}

\begin{lemma}
\label{il07}With the above notations, any CM-subfield $L$ of $H_{\mathfrak{m}%
{}}(K{})$ containing $K$ is fixed by the subgroup $2\cdot C_{\mathfrak{m}{}%
}^{-}(K)$ of $C_{\mathfrak{m}{}}(K)$, i.e., $2\cdot C_{\mathfrak{m}{}}%
^{-}(K)\hookrightarrow\Gal(H_{\mathfrak{m}{}}(K)/L)$.
\end{lemma}

\begin{proof}
We shall need the general fact:%
\begin{quote}%
(*). Let $L/F$ be a finite Galois extension of number fields, and let $\tau$
be a homomorphism $L\rightarrow$ $\mathbb{Q}^{\text{al}}$. For any prime ideal
$\mathfrak{P}{}$ of $L$ unramified in $L/F$,
\[
(\tau\mathfrak{P}{},\tau L/\tau F)=\tau\circ(\mathfrak{\mathfrak{P}{}%
},L/F)\circ\tau^{-1}%
\]
(equality of Artin symbols). In particular, if $L/F$ is abelian, $\tau L=L$,
and $\tau F=F$, then for any prime ideal $\mathfrak{p}{}$ of $F$ unramified in
$L$, $(\tau\mathfrak{\mathfrak{p}{}},L/F)=\tau\circ(\mathfrak{p},L/F)\circ
\tau^{-1}$.%
\end{quote}%
Because $\iota_{K}\mathfrak{m}{}=\mathfrak{m}{}$, $H_{\mathfrak{m}{}}(K)$ is
stable under $\iota$, and hence is Galois over $K_{+}$. The Galois closure of
$L$ over $K_{+}$ will again be CM, and so we may assume $L$ to be Galois over
$K_{+}$. Let $H$ be the subgroup of $C_{\mathfrak{m}{}}(K)$ fixing $L$, and
consider the diagram%
\[
\begin{CD} 0 @>>> C_{\mathfrak{m}}(K)@>>> \Gal(H_{\mathfrak{m}}(K)/K_{+}) @>>> \Gal(K/K_{+}) @>>> 0 \\ @.@VVV@VVV@| \\ 0 @>>> C_{\mathfrak{m}}(K)/H @>>> \Gal(L/K_{+})@>>> \Gal(K/K_{+}) @>>>0 \end{CD}
\]
in which we have used the reciprocity map to replace $\Gal(H_{\mathfrak{m}{}%
}(K)/K)$ with $C_{\mathfrak{m}{}}(K)$. According to (*), the action of
$\Gal(K/K_{+})=\langle\iota_{K}\rangle$ on $C_{\mathfrak{m}{}}(K)$ defined by
the upper sequence is the natural action. Because $L$ is a CM-field, the
action of $\iota_{K}$ on $C_{\mathfrak{m}{}}(K)/H$ defined by the lower
sequence is trivial. The composite of the maps
\[
\begin{CD}
C_{\mathfrak{m}}(K)@>{\Nm_{K/K_{+}}}>>C_{\mathfrak{m}{}%
_{+}}(K_{+})@>>> C_{\mathfrak{m}{}}(K)
\end{CD}
\]
is $c\mapsto c\cdot\iota c$, and so $\iota_{K}$ acts as $-1$ on
$C_{\mathfrak{m}{}}^{-}(K)$. Thus $H\supset2C_{\mathfrak{m}{}}^{-}(K)$.
\end{proof}

For a finite abelian group $A$, let $A(\mathrm{odd})$ denote the subgroup of
elements of odd order in $A$.

\begin{lemma}
\label{il08}For any CM-field $L$ containing $K$, the norm map
\[
C^{-}(L)(\mathrm{odd})\rightarrow C^{-}(K)(\mathrm{odd})
\]
is surjective.
\end{lemma}

\begin{proof}
The norm limitation theorem (\cite{tate1967}, p202, paragraph 3) allows us to
assume that $L$ is abelian over $K$. Therefore, we may suppose $L\subset
H_{\mathfrak{m}{}}(K)$ for some modulus $\mathfrak{m}{}$ as above. Moreover,
we may suppose that $L$ is the largest CM-subfield of $H_{\mathfrak{m}{}}(K)$.
The field $H_{\mathfrak{m}{}_{+}}(K_{+})\cdot K$ is CM, and so is contained in
$L$. The field $H_{\mathfrak{m}{}_{+}}(L_{+})\cdot L$ is also CM, and so
\begin{equation}
(H_{\mathfrak{m}{}_{+}}(L_{+})\cdot L)\cap H_{\mathfrak{m}{}}(K)=L.\label{e9}%
\end{equation}
Consider the diagram (cf. (\ref{1d}))
\[
\begin{CD} C_{\mathfrak{m}}^{-}(L) @>{\cong}>> \Gal(H_{\mathfrak{m}}(L)/H_{\mathfrak{m}_{+}}(L_{+})\cdot L) \\ @VV{\Nm_{L/K}}V @VV{\sigma \mapsto \sigma |H_{\mathfrak{m}}(K)}V \\ C_{\mathfrak{m}}^{-}(K) @>{\cong}>> \Gal(H_{\mathfrak{m}}(K)/H_{\mathfrak{m}_{+}}(K_{+})\cdot K) \\ \end{CD}
\]
According to (\ref{e9}), the image of $\sigma\mapsto\sigma|H_{\mathfrak{m}{}%
}(K)$ is $\Gal(H_{\mathfrak{m}{}}(K)/L)$, and so by (\ref{il07}) the image of
$\Nm_{L/K}$ contains $2\cdot C_{\mathfrak{m}{}}^{-}(K)$. Therefore,
$\Nm_{L/K}\colon C_{\mathfrak{m}{}}^{-}(L)($odd$)\rightarrow C_{\mathfrak{m}%
{}}^{-}(K)($odd$)$ is surjective, which implies the similar statement without
the $\mathfrak{m}{}$'s.
\end{proof}

\begin{proof}
[\textsc{Proof of Proposition \ref{il06}}]As $\plim{}C(K_{+})=0$
(\ref{il05m}), the inverse system of exact sequences (\ref{e3}) gives an
isomorphism%
\[
\plim{}C^{-}(K)\cong\plim{}C(K)\text{.}%
\]
Because of (\ref{il08}), for every CM-field $K_{0}$, the projection%
\[
\plim{}C^{-}(K)(\text{odd})\rightarrow C^{-}(K_{0})(\text{odd})
\]
is surjective. An irregular prime $l$ divides the order of $C^{-}(\mathbb{Q}%
{}[\zeta_{l}])$ (\cite{washington1997}, 5.16, p62). Therefore the order of the
profinite group $\plim{}C^{-}(K)($odd$)$ is divisible by every irregular
prime. Since there are infinitely many irregular primes (ibid. 5.17), this
implies that $\plim{}C^{-}(K)($odd$)$ is infinite and hence
uncountable.\footnote{Let $h_{K}^{-}$ be the order of $C^{-}(K)$ (the relative
class number). According to the Brauer-Siegel theorem, as $K$ ranges over CM
number fields of a given degree, $h_{K}^{-}$ is asymptotic to $\log
(\sqrt{d_{K}/d_{K_{+}}})$ and goes to infinity with $d_{K}$ (the discriminant
of $K$). The Generalized Riemann Hypothesis implies a similar statement for
the exponents of the class groups (\cite{louboutin2003}).}
\end{proof}

\subsection{The inverse system $(K^{\times}/U(K))$}

\begin{lemma}
\label{il02}For any set $\mathcal{F}{}$ satisfying (\ref{il01}a,b), there is
an exact sequence%
\[
0\rightarrow\plim{}_{\mathcal{F}{}}F^{\times}/U(F)\rightarrow
\plim{}_{\mathcal{F}{}}I(F)\rightarrow\plim{}_{\mathcal{F}{}}C(F)\rightarrow
\plim{}_{\mathcal{F}{}}^{1}F^{\times}/U(F)\rightarrow\plim{}_{\mathcal{F}{}%
}^{1}I(F)\rightarrow0\text{.}%
\]

\end{lemma}

\begin{proof}
Because the groups $C(F)$ are finite, $\plim{}^{1}C(F)=0$ (\ref{il01b}), and
so this is the exact sequence (\ref{e1}) attached to the inverse system of
short exact sequences%
\[
0\rightarrow F^{\times}/U(F)\rightarrow I(F)\rightarrow C(F)\rightarrow
0\text{.}%
\]

\end{proof}

\begin{lemma}
\label{il03}Let $\mathcal{F}{}$ be the set of totally real fields in
$\mathbb{Q}{}^{\text{al}}$ of finite degree over $\mathbb{Q}{}$. Then%
\begin{align}
\plim{}_{\mathcal{F}{}}F^{\times} &  =0\label{e26}\\
\plim{}_{\mathcal{F}{}}\mathbb{I}{}(F) &  \cong\plim{}_{\mathcal{F}{}%
}(F\otimes_{\mathbb{Q}{}}\mathbb{R}{})^{\times}\label{e27}\\
\plim{}_{\mathcal{F}{}}F^{\times}/U(F) &  =0\text{,}\label{e28}\\
\plim{}_{\mathcal{F}{}}^{1}F^{\times}/U(F) &  \cong\plim{}_{\mathcal{F}{}}%
^{1}I(F).\label{e29}%
\end{align}
Moreover, $\plim{}_{\mathcal{F}{}}^{1}F^{\times}/U(F)$ is uncountable.
\end{lemma}

\begin{proof}
The equality (\ref{e27}) follows from the discussion (\ref{il05}) and the fact
that a nonarchimedean local field has no universal norms. Equality (\ref{e27})
implies (\ref{e26}) because a global norm is a local norm.

An ideal $\mathfrak{a}{}=\mathfrak{p}{}^{m}\cdots\in I(F)$, $m\neq0$, will not
be in the image of
\[
\Nm_{L/F}\colon I(L)\rightarrow I(F)
\]
once $L$ is so large that none of the residue class degrees of the primes
above $\mathfrak{p}{}$ divides $m$. It follows that $\plim{}I(F)=0$ and that
$(I(F))_{F\in\mathcal{F}{}}$ fails (ML). Thus, (\ref{il02}) proves (\ref{e28})
and that there is an exact sequence%
\[
0\rightarrow\plim{}C(F)\rightarrow\plim{}^{1}F^{\times}/U(F)\rightarrow
\plim{}^{1}I(F)\rightarrow0\text{.}%
\]
But $\varprojlim C(F)=0$ (\ref{il05m}) and so $\varprojlim^{1}F^{\times
}/U(F)\cong\varprojlim^{1}I(F)$, which is uncountable (\ref{il01b}).
\end{proof}

\begin{lemma}
\label{il04}Let $\mathcal{F}{}$ be the set of CM-subfields of $\mathbb{Q}%
^{\text{al}}$ of finite degree over $\mathbb{Q}{}$. Then%
\begin{align}
\plim{}_{\mathcal{F}{}}K^{\times} &  =0\label{e30}\\
\plim{}_{\mathcal{F}{}}\mathbb{I}{}(K) &  \cong\plim{}_{\mathcal{F}{}%
}(K\otimes_{\mathbb{Q}{}}\mathbb{R})^{\times}\label{e31}\\
\plim{}_{\mathcal{F}{}}K^{\times}/U(K) &  =0\label{e32}%
\end{align}
and there is an exact sequence%
\begin{equation}
0\rightarrow\plim{}_{\mathcal{F}{}}{}C(K)\rightarrow\plim{}_{\mathcal{F}{}%
}^{1}K^{\times}/U(K)\rightarrow\plim{}_{\mathcal{F}{}}^{1}I(K)\rightarrow
0\text{.}\label{e33}%
\end{equation}

\end{lemma}

\begin{proof}
The proofs of (\ref{e30}), (\ref{e31}), and (\ref{e32}) are similar to those
of the corresponding equalities in (\ref{il03}). In particular,
$\plim{}I(K)=0$, and so (\ref{e33}) follows from \ref{il02}.
\end{proof}

\subsection{The inverse system $(K^{\times}/K_{+}^{\times})$}

Let $\mathcal{F}{}$ be the set of CM fields in $\mathbb{Q}{}^{\text{al}}$ of
finite degree over $\mathbb{Q}{}$.

\begin{lemma}
\label{il09}The quotient map%
\[
\plim{}_{\mathcal{F}{}}^{1}K^{\times}/K_{+}^{\times}\rightarrow
\plim{}_{\mathcal{F}{}}^{1}(K^{\times}/U(K)\cdot K_{+}^{\times})
\]
is an isomorphism.
\end{lemma}

\begin{proof}
The homomorphism $U(K_{+})\rightarrow U(K)$ has finite cokernel, and so
\[
\plim{}^{1}U(K_{+})\rightarrow\plim{}^{1}U(K)
\]
is surjective. Now the statement follows from the diagram%
\[
\begin{CD} \plim{}^{1}U(K_+) @>>> \plim{}^{1}K_+^{\times }
@>>> \plim{}^{1}K_+^{\times }/U(K_+) @>>> 0 \\
@VVV @VVV @VVV \\
\plim{}^{1}U(K) @>>> \plim{}^{1}K^{\times } @>>>
\plim{}^{1}K^{\times }/U(K) @>>>0\\
@VVV@VVV@VVV\\
0@>>>\plim{}^{1}K^{\times}/K_+^{\times}@>>>
\plim{}^{1}K^{\times}/U(K)\cdot K_+^{\times}\\
@.@VVV@VVV\\
@.0@.0@..\end{CD}
\]

\end{proof}

\begin{lemma}
\label{il09m}The quotient map%
\[
\plim{}_{\mathcal{F}{}}^{1}\mathbb{I}{}(K)/\mathbb{I}{}(K_{+})\rightarrow
\plim{}_{\mathcal{F}{}}^{1}I(K)/I(K_{+})
\]
is an isomorphism.
\end{lemma}

\begin{proof}
For a number field $L$, let $V(L)$ denote the kernel of $\mathbb{I}%
{}(L)\rightarrow I(L)$. Let $r=[K_{+}\colon\mathbb{Q}{}]$. Then $V(K_{+}%
)\approx(\mathbb{R}{}^{\times})^{r}\times\left\{  \text{\textrm{compact}%
}\right\}  \,$and $V(K)\approx$ $(\mathbb{\mathbb{C}{}}^{\times})^{r}%
\times\left\{  \text{\textrm{compact}}\right\}  $, and so $V(K)/V(K_{+})$ is
compact. Therefore, $\varprojlim^{1}V(K)/V(K_{+})=1$ (\ref{il01b}), and the
limit of the exact sequences%
\[
0\rightarrow V(K)/V(K_{+})\rightarrow\mathbb{I}{}(K)/\mathbb{I}{}%
(K_{+})\rightarrow I(K)/I(K_{+})\rightarrow0
\]
gives the required isomorphism.
\end{proof}

\begin{proposition}
\label{il10}There is an exact sequence%
\[
0\rightarrow\plim{}_{\mathcal{F}{}}C(K)\rightarrow\plim{}_{\mathcal{F}{}}%
^{1}K^{\times}/K_{+}^{\times}\rightarrow\plim{}_{\mathcal{F}{}}^{1}%
\mathbb{I}{}(K)/\mathbb{I}{}(K_{+})\rightarrow0\text{.}%
\]

\end{proposition}

\begin{proof}
Lemmas \ref{il03} and \ref{il04} show that the rows in%
\[
\begin{CD} @. 0 @>>> \plim{}_{\mathcal{F}}^{1}K_{+}^{\times }/U(K_{+}) @>>>\plim{}_{\mathcal{F}{}}^{1}I(K_{+}) @>>> 0 \\ @.@.@VVV@VVV@. \\ 0 @>>> \plim{}_{\mathcal{F}}C(K) @>>> \plim{}_{\mathcal{F}}^{1}K^{\times }/U(K) @>>> \plim{}_{\mathcal{F}}^{1}I(K) @>>> 0\\ \end{CD}
\]
are exact, from which we deduce an exact sequence%
\[
0\rightarrow\plim{}_{\mathcal{F}{}}{}C(K)\rightarrow\plim{}_{\mathcal{F}{}%
}^{1}(K^{\times}/U(K)\cdot K_{+}^{\times})\rightarrow\plim{}_{\mathcal{F}{}%
}^{1}I(K)/I(K_{+})\rightarrow0\text{.}%
\]
Now (\ref{il09}) and (\ref{il09m}) allow us to replace the last two terms with
$\plim{}_{\mathcal{F}{}}^{1}K^{\times}/K_{+}^{\times}$ and
$\plim{}_{\mathcal{F}{}}^{1}\mathbb{I}{}(K)/\mathbb{I}{}(K_{+})$.
\end{proof}

\subsection{The inverse system $(K(w)^{\times}/K(w)_{+}^{\times})$}

Fix a finite prime $w$ of $\mathbb{Q}{}^{\text{al}}$, and write $w_{L}$ (or
just $w$) for the prime it defines on a subfield $L$ of $\mathbb{Q}%
^{\text{al}}$. For an $L$ that is Galois over $\mathbb{Q}{}$, let $D(w_{L})$
denote the decomposition group, and let $L(w)=L^{D(w_{L})}$. For a system
$\mathcal{F}{}$ satisfying (\ref{il01}) and whose members are Galois over
$\mathbb{Q}{}$, we define $(L(w)^{\times})_{L\in\mathcal{F}{}}$ to be the
inverse system with transition maps
\[
\left(  \Nm_{L^{\prime}(w)/L(w)}\right)  ^{[L_{w}^{\prime}\colon L_{w}]}\colon
L^{\prime}(w)^{\times}\rightarrow L(w)^{\times}\text{.}%
\]
We define the inverse systems $\left(  C(L(w))\right)  _{F\in\mathcal{F}{}}$,
$(I(L(w))_{F\in\mathcal{F}{}}$, etc., similarly.

\begin{proposition}
\label{il11}Let $\mathcal{F}{}$ be the set of totally real fields in
$\mathbb{Q}{}^{\text{al}}$ of finite degree and Galois over $\mathbb{Q}{}$.
Then%
\[
\plim{}_{\mathcal{F}{}}(C(F(w))=0
\]
and the map%
\[
\plim{}_{\mathcal{F}{}}^{1}F(w)^{\times}/U(F(w))\rightarrow
\plim{}_{\mathcal{F}{}}^{1}I(F(w))
\]
is an isomorphism.
\end{proposition}

\begin{proof}
Let $m$ be the exponent of $C(F(w))$. The $F^{\prime}$ in $\mathcal{F}$ such
that $F^{\prime}\supset F$ and $m|[F_{w}^{\prime}\colon F_{w}]$ form a cofinal
subset of $\mathcal{F}{}$.\footnote{Let $f_{0}$ be a monic polynomial in
$\mathbb{Q}{}[X]$ whose splitting field is $F$; let $L$ be a finite extension
of $F_{w}$ such that $h|[L\colon F_{w}]$, and write $L\approx\mathbb{Q}{}%
_{p}[X]/(g(X))$ with $g$ monic; choose a monic $f\in\mathbb{Q}{}[X]$ such that
$f$ is close $w$-adically to $g$ and really close to a monic polynomial that
splits over $\mathbb{R}{}$; then, by Krasner's lemma, the splitting field
$F_{0}^{\prime}$ of $f_{0}\cdot f$ will have the required property, as will
any field in $\mathcal{F}{}$ containg $F_{0}^{\prime}$.} For such an
$F^{\prime}$, the map
\[
\left(  \Nm_{F^{\prime}(w)/F(w)}\right)  ^{[F_{w}^{\prime}\colon F_{w}]}\colon
C(F^{\prime}(w))\rightarrow C(F(w))
\]
is zero, and so $\plim{}_{\mathcal{F}{}}(C(F(w))=0$ by the criterion
(\ref{il01g}). The second statement now follows from (\ref{il02}).
\end{proof}

\begin{proposition}
\label{il12}Let $\mathcal{F}{}$ be the set of CM fields in $\mathbb{Q}%
^{\text{al}}$ of finite degree and Galois over $\mathbb{Q}{}$. Then%
\begin{equation}
\plim{}_{\mathcal{F}{}}{}C(K(w))=0\label{e20}%
\end{equation}
and%
\begin{align}
\plim{}_{\mathcal{F}{}}^{1}K(w)^{\times}/U(K(w)) &  \cong\plim{}_{\mathcal{F}%
{}}^{1}I(K(w))\newline\label{e21}\\
\plim{}_{\mathcal{F}{}}^{1}K(w)^{\times}/K(w)_{+}^{\times} &  \cong
\plim{}_{\mathcal{F}{}}^{1}K(w)^{\times}/U(K(w))\cdot K(w)_{+}^{\times
}\newline\label{e22}\\
\plim{}_{\mathcal{F}{}}^{1}\mathbb{I}{}(K(w))/\mathbb{I}{}(K(w)_{+}) &
\cong\plim{}_{\mathcal{F}{}}^{1}I(K(w))/I(K(w)_{+})\newline\label{e23}\\
\plim{}_{\mathcal{F}{}}^{1}K(w)^{\times}/K(w)_{+}^{\times} &  \cong
\plim{}_{\mathcal{F}{}}^{1}\mathbb{I}{}(K(w))/\mathbb{I}{}(K(w)_{+}%
).\label{e24}%
\end{align}

\end{proposition}

\begin{proof}
The proof of (\ref{e20}) is similar to that of (\ref{il11}) (cf. \ref{il05}b),
and (\ref{e21}) then follows from (\ref{il02}). The proofs of (\ref{e21}) and
(\ref{e22}) are similar to those of (\ref{il09}) and (\ref{il09m}). For
(\ref{e24}), the horizontal arrows in the commutative diagram%
\[
\begin{CD} \plim{}_{\mathcal{F}}^{1}K(w)^{\times }_+/U(K(w)_+)
@>{\cong}>> \plim{}_{\mathcal{F}}^{1}I(K(w)_+) \\
@VVV@VVV\\
\plim{}_{\mathcal{F}}^{1}K(w)^{\times }/U(K(w)) @>{\cong}>> \plim{}_{\mathcal{F}}^{1}I(K(w)) .
\end{CD}
\]
are isomorphisms by Proposition \ref{il11} and (\ref{e20}). Thus,
\[
\plim{}_{\mathcal{F}{}}^{1}(K(w)^{\times}/U(K(w))\cdot K(w)_{+}^{\times}%
)\cong\plim{}_{\mathcal{F}{}}^{1}I(K(w))/I(K(w)_{+}),
\]
and (\ref{il12}b) and (\ref{il12}c) allow us to replace the terms with
$\plim{}_{\mathcal{F}{}}^{1}K(w)^{\times}/K(w)_{+}^{\times}$ and
$\plim{}_{\mathcal{F}{}}^{1}\mathbb{I}{}(K(w))/\mathbb{I}{}(K(w)_{+})$.
\end{proof}

\subsection{Conclusions}

\begin{proposition}
\label{il13}Let $\mathcal{F}{}$ be the set of CM fields in $\mathbb{Q}%
^{\text{al}}$ of finite degree and Galois over $\mathbb{Q}{}$. In the
diagram,
\[
\begin{CD}
\varprojlim_{\mathcal{F}}^{1}K(w)^{\times}/K(w)_{+}^{\times}@>>>
\varprojlim_{\mathcal{F}}^{1}\mathbb{I}(K(w))/\mathbb{I}%
(K(w)_{+}\\
@VVcV@VVV\\
\varprojlim_{\mathcal{F}{}}^{1}K^{\times}/K_{+}^{\times}
@>d>>\varprojlim_{\mathcal{F}{}}^{1}\mathbb{I}{}(K)/\mathbb{I}{}(K^{+}),
\end{CD}
\]
$\Ker(d\circ c)=0$ and $\Ker(d)\cong\varprojlim C(K)$.
\end{proposition}

\begin{proof}
The second statement is proved in (\ref{il10}). Since the top horizontal arrow
is an isomorphism (\ref{e24}), for the first statement, it suffices to show
that right vertical map is injective. Because of (\ref{il09m}) and
(\ref{e23}), this is equivalent to showing that
\[
\varprojlim{}^{1}I(K(w))/I(K(w)_{+})\rightarrow\varprojlim{}^{1}I(K)/I(K_{+})
\]
is injective. From the exact sequence (\ref{e1}) attached to%
\[
0\rightarrow I(K(w))/I(K(w)_{+}\rightarrow I(K)/I(K_{+})\rightarrow
I(K)/\left(  I(K_{+})+I(K(w))\right)  \rightarrow0
\]
we see that it suffices to show that
\[
\varprojlim I(K)/\left(  I(K_{+})+I(K(w))\right)  =0.
\]
Let $a\in I(K)$ represent a nonzero element of $I(K)/\left(  I(K_{+}%
)+I(K(w))\right)  $. For some $m$, $a\notin I(K_{+})+I(K(w))+mI(K)$. There
exists a $K^{\prime}\in\mathcal{F}{}$ containing $K$ and such that, for every
prime $v^{\prime}$ of $K^{\prime}$ lying over a prime $v$ of $K$ for which the
$v$-component of $a$ is nonzero, $m$ divides the residue class degree
$f(v^{\prime}/v)$. For such a $K^{\prime}$, $a\notin I(K_{+})+I(K(w))+\Nm
I(K^{\prime})$. Thus, $\varprojlim I(K)/I(K_{+})\cdot I(K(w))=0$ by the
criterion (\ref{il01g}).
\end{proof}

\newpage

\section{The cohomology of protori}

Throughout this section, $k$ is a field of characteristic zero, $k^{\text{al}%
}$ is an algebraic closure of $k$, and $\Gamma=\Gal(k^{\text{al}}/k)$.

\subsection{Review of affine group schemes}

\begin{plain}
\label{ra1}Every affine group scheme over $k$ is the inverse limit of a strict
inverse system of algebraic groups (\cite{demazure1970}, III, \S 3, 7.5, p355).
\end{plain}

\begin{plain}
\label{ra2}An affine group scheme $G$ over $k$ defines a sheaf $\tilde
{G}\colon U\mapsto G(U)$ on $(\Spec k)_{\text{fpqc }}$, and the functor
$G\mapsto\tilde{G}$ is fully faithful (ibid. III, \S 1, 3.3, p297). When $N$
is an affine normal subgroup scheme of $G$, the quotient sheaf $\tilde
{G}/\tilde{N}$ is in the essential image of the functor (ibid. III, \S 3, 7.2, p353).
\end{plain}

\begin{plain}
\label{ra2m}Let $(G_{\alpha})_{\alpha\in I}$ be an inverse system of affine
group schemes over $k$. Then $G=\varprojlim G_{\alpha}$ is an affine group
scheme, and it is the inverse limit of $(G_{\alpha})_{\alpha\in I}$ in the
category of $k$-schemes, i.e., $G(S)=\varprojlim G_{\alpha}(S)$ for all
$k$-schemes $S$ (ibid. III, \S 3, 7.5, p355). It follows\footnote{Recall that
to form an inverse limit in the category of sheaves, form the inverse limit in
the category of presheaves (this is the obvious object), and then the
resulting presheaf is a sheaf and is the inverse limit in the category of
sheaves.} that $\tilde{G}=\varprojlim\tilde{G}_{\alpha}$.
\end{plain}

\begin{plain}
\label{ra3}The category of commutative affine group schemes over $k$ is
abelian. A sequence%
\[
0\rightarrow A\overset{a}{\rightarrow}B\overset{b}{\rightarrow}C\rightarrow0
\]
is exact if and only if $A\rightarrow B$ is a closed immersion, $B\rightarrow
C$ is fully faithful, and $a\colon A\rightarrow B$ is a kernel of $b$ (ibid.
III, \S 3, 7.4, p355). The category is pro-artinian, and so inverse limits of
exact sequences are exact (ibid. V, \S 2, 2, pp563--5).
\end{plain}

\begin{plain}
\label{ra4}The functor $G\mapsto\tilde{G}$ from the category of commutative
affine group schemes over $k$ to the category of sheaves of commutative groups
on $(\Spec k)_{\text{fpqc }}$ is exact (left exactness is obvious, and right
exactness follows from \ref{ra2}).
\end{plain}

\begin{plain}
\label{ra5}We say that an affine group scheme $G$ is \emph{separable}%
\footnote{Compare: a topological space is said to be separable if it has a
countable dense subset; a profinite group $G$ is separable if and only if the
set of its open subgroups is countable, in which case $G$ is the limit of a
strict inverse system of finite groups indexed by $(\mathbb{N}{},\leq)$%
.}\emph{ }if the set of affine normal subgroups $N$ of $G$ for which $G/N$
algebraic is countable. Such a $G$ is the inverse limit of a strict inverse
system of algebraic groups indexed by $(\mathbb{N}{},\leq)$ (apply \ref{ra1}).
\end{plain}

\subsection{Continuous cohomology}

Let $T$ be a separable protorus over $k$. Then%
\[
T(k^{\text{al}})=\varprojlim Q(k^{\text{al}})
\]
where $Q$ runs over the algebraic quotients of $T$. Endow each group
$Q(k^{\text{al}})$ with the discrete topology, and endow $T(k^{\text{al}})$
with the inverse limit topology. Define $H_{\text{cts}}^{r}(k,T)$ to be the
$r$th cohomology group of the complex%
\[
\cdots\rightarrow C^{r}(\Gamma,T)\xr{d^r}C^{r+1}(\Gamma,T)\rightarrow\cdots
\]
where $C^{r}(\Gamma,T)$ is the group of continuous maps $\Gamma^{r}\rightarrow
T(k^{\text{al}})$ and $d^{r}$ is the usual boundary map (see, for example,
\cite{tate1976}, \S 2). Note that, if $T=\varprojlim_{n}T_{n}$, then
$C^{r}(\Gamma,T)=\varprojlim_{n}C^{r}(\Gamma,T_{n})$.

When $T$ is a torus, $H_{\text{cts}}^{r}(k,T)$ is the usual Galois cohomology
group%
\[
H_{\text{cts}}^{r}(k,T)=H^{r}(\Gamma,T(k^{\text{al}}))=\dlim H^{r}%
(\Gamma_{K/k},T(K))
\]
(limit over the finite Galois extensions $K$ of $k$ contained in
$k^{\text{al}}$; $\Gamma_{K/k}=\Gal(K/k)$). In this case, we usually omit the
subscript \textquotedblleft cts\textquotedblright.

\begin{lemma}
\label{cc1}Let
\begin{equation}
0\rightarrow T^{\prime}\rightarrow T\rightarrow T^{\prime\prime}\rightarrow0
\label{e11}%
\end{equation}
be the limit of a short exact sequence of countable strict inverse systems of
tori. Then the sequence of complexes%
\begin{equation}
0\rightarrow C^{\bullet}(\Gamma,T^{\prime})\rightarrow C^{\bullet}%
(\Gamma,T)\rightarrow C^{\bullet}(\Gamma,T^{\prime\prime})\rightarrow0
\label{e13}%
\end{equation}
is exact, and so gives rise to a long exact sequence%
\[
\cdots\rightarrow H_{\text{cts}}^{r}(k,T^{\prime})\rightarrow H_{\text{cts}%
}^{r}(k,T)\rightarrow H_{\text{cts}}^{r}(k,T^{\prime\prime})\rightarrow
H_{\text{cts}}^{r+1}(k,T^{\prime})\rightarrow\cdots.
\]

\end{lemma}

\begin{proof}
By assumption, (\ref{e11}) is the inverse limit of a system%
\[
\begin{CD}
@.\vdots@.\vdots@.\vdots\\
0@>>>T_{n+1}^{\prime}@>>>T_{n+1}@>>>T_{n+1}^{\prime\prime}@>>>0\\
@.@VV{\mathrm{onto}}V@VV{\mathrm{onto}}V@VV{\mathrm{onto}}V@.\\
0@>>>T_{n}^{\prime}@>>>T_{n}@>>>T_{n}^{\prime\prime}@>>>0\\
@.\vdots@.\vdots@.\vdots\\
\end{CD}
\]
Because the transition maps $T_{n+1}^{\prime}(k^{\text{al}})\rightarrow
T_{n}^{\prime}(k^{\text{al}})$ are surjective, the limit of the short exact
sequences%
\[
0\rightarrow T_{n}^{\prime}(k^{\text{al}})\rightarrow T_{n}(k^{\text{al}%
})\rightarrow T_{n}^{\prime\prime}(k^{\text{al}})\rightarrow0
\]
is an exact sequence%
\[
0\rightarrow T^{\prime}(k^{\text{al}})\rightarrow T(k^{\text{al}})\rightarrow
T^{\prime\prime}(k^{\text{al}})\rightarrow0
\]
(\ref{il01b}(b)). To show that (\ref{e13}) is exact, it suffices to show that

\begin{enumerate}
\item the topology on $T^{\prime}(k^{\text{al}})$ is induced from that on
$T(k^{\text{al}})$;

\item the map $T(k^{\text{al}})\rightarrow T^{\prime\prime}(k^{\text{al}})$
has a continuous section (not necessarily a homomorphism).
\end{enumerate}

\noindent(see \cite{tate1976}, p259).

\noindent The topology on an inverse limit is that inherited by it as a subset
of the product. Thus (a) follows from the similar statement for products
(\cite{bourbaki1989}, I 4.1, Corollary to Proposition 3, p46). As all the
groups in the diagram%
\[
\begin{CD}
T_{n}(k^{\text{al}}) @>>> T_{n}^{\prime\prime}(k^{\text{al}})\\
@AAA@AAA\\
T_{n+1}(k^{\text{al}}) @>>> T_{n+1}^{\prime\prime}(k^{\text{al}})
\end{CD}
\]
are discrete and all the maps surjective, it is possible to successively
choose compatible sections to the maps $T_{n}(k^{\text{al}})\rightarrow
T_{n}^{\prime\prime}(k^{\text{al}})$. Their limit is the section required for (b).
\end{proof}

\begin{proposition}
\label{cc2}Let $T$ be a separable protorus over $k$, and write it as the limit
$T=\varprojlim_{n}T_{n}$ of a strict inverse system of tori. For each $r\geq
0$, there is an exact sequences%
\[
0\rightarrow\plim{}^{1}H^{r-1}(k,T_{n})\rightarrow H_{\text{cts}}%
^{r}(k,T)\rightarrow\varprojlim H^{r}(k,T_{n})\rightarrow0\text{.}%
\]

\end{proposition}

\begin{proof}
The map $T_{n+1}(k^{\text{al}})\rightarrow T_{n}(k^{\text{al}})$ is
surjective, and admits a continuous section because $T_{n}(k^{\text{al}})$ is
discrete. Hence $C^{r}(\Gamma,T_{n+1})\rightarrow C^{r}(\Gamma,T_{n})$ is
surjective. Thus, by (\ref{il01b}),%
\[
\plim{}^{1}C^{r}(\Gamma,T_{n})=0,
\]
and so we can apply the next lemma to the inverse system of complexes
$(C^{\bullet}(\Gamma,T_{n}))_{n\in\mathbb{N}{}}$.
\end{proof}

\begin{lemma}
\label{cc3}Let $(C_{n}^{\bullet})_{n\in\mathbb{N}{}}$ be an inverse system of
complexes of abelian groups such that $\varprojlim_{n}^{1}C_{n}^{r}=0$ for all
$r$. Then there is a canonical exact sequence%
\begin{equation}
0\rightarrow\plim{}_{n}^{1}H^{r-1}(C_{n}^{\bullet})\rightarrow H^{r}%
(\varprojlim C_{n}^{\bullet})\rightarrow\varprojlim H^{r}(C_{n}^{\bullet
})\rightarrow0. \label{e14}%
\end{equation}

\end{lemma}

\begin{proof}
This is a standard result. [The condition $\varprojlim_{n}^{1}C_{n}^{r}=0$
implies that the sequence of complexes%
\[
\begin{CD} @.\vdots@.\vdots@.\vdots@.\\
@.@VVV@VVV@VVV@.\\
0@>>>\plim_{n}C_{n}^{r}@>>>
\textstyle\prod\nolimits_{n}C^{r}_n@>1-u>>
\textstyle\prod\nolimits_{n}C^{r}_n@>>>0\\
@.@VVV@VVV@VVV@.\\
0@>>>\plim_{n}C_{n}^{r+1} @>>>
\textstyle\prod\nolimits_{n}C^{r+1}_n@>1-u>>
\textstyle\prod\nolimits_{n}C^{r+1}_n@>>>0\\
@.@VVV@VVV@VVV@.\\ @.\vdots@.\vdots@.\vdots@.\\ \end{CD}
\]
is exact (the maps $1-u$ are as in (\ref{e10})). The associated long exact
sequence is
\[
\cdots\xr{1-u}%
{\textstyle\prod\nolimits_{n}}
H^{r-1}(C_{n}^{\bullet})\rightarrow H^{r}(\varprojlim C_{n}^{\bullet
})\rightarrow%
{\textstyle\prod\nolimits_{n}}
H^{r}(C_{n}^{\bullet})\xr{1-u}%
{\textstyle\prod\nolimits_{n}}
H^{r}(C_{n}^{\bullet})\rightarrow\cdots,
\]
which gives (\ref{e14}).]
\end{proof}

\begin{lemma}
\label{cc4}For a separable protorus $T$, $H_{\text{fpqc}}^{1}(k,T)\cong
H_{\text{cts}}^{1}(k,T)$.
\end{lemma}

\begin{proof}
The group $H_{\text{fpqc}}^{1}(k,T)$ is canonically isomorphic to the group of
isomorphism classes of torsors for $T$ (\cite{giraud1971}, III 3.5.4, p169).
Let $P$ be such a torsor. Because $T$ is separable, it is the limit of a
strict inverse system $(T_{n})_{n\in\mathbb{N}{}}$, and correspondingly $P$ is
the inverse limit $P=\varprojlim P_{n}$ of the strict inverse system
$(P_{n})_{n\in\mathbb{N}{}}=\left(  P\wedge^{T}T_{n}\right)  _{n\in
\mathbb{N}{}}$. The set $P(k^{\text{al}})=\varprojlim P_{n}(k^{\text{al}})$ is
nonempty. Choose a $p\in P(k^{\text{al}})$, and for $\sigma\in\Gamma$ write
$\sigma p=p\cdot a_{\sigma}(P)$ with $a_{\sigma}(P)\in T(k^{\text{al}})$. Then
$\sigma\mapsto a_{\sigma}(P)$ is a cocycle, which is continuous because its
projection to each $T_{n}(k^{\text{al}})$ is continuous, and its class in
$H_{\text{cts}}^{1}(k,T)$ depends only on the isomorphism class of $P$. The
map
\[
\lbrack P]\mapsto\lbrack a_{\sigma}(P)]\colon H_{\text{fpqc}}^{1}%
(k,T)\rightarrow H_{\text{cts}}^{1}(k,T)
\]
is easily seen to be injective, and the flat descent theorems show that it is
surjective. (See \cite{milne2002c}, 1.20, for a slightly different proof.)
\end{proof}

\begin{example}
\label{cc5}Define the universal covering $\tilde{T}$ of a torus $T$ to be the
inverse limit of the inverse system indexed by $\mathbb{N}{}^{\times}$%
\[
\cdots\leftarrow\underset{n}{T}\,\overset{m}{\leftarrow}\,\underset{mn}%
{T}\,\leftarrow\cdots.
\]
Then, with the notations of (\ref{il01p}), there is an exact sequence%
\[
0\rightarrow\plim{}^{1}(H^{r-1}(k,T),m)\rightarrow H^{r}(k,\tilde
{T})\rightarrow\varprojlim(H^{r}(k,T),m)\rightarrow0.
\]
For any torus $T$ over $\mathbb{Q}{}$, $\tilde{T}(\mathbb{Q}{})=0$ (e.g.,
\cite{milne1994a}, 3.16).\footnote{An element of $\widetilde{T}(\mathbb{Q})$
is a family $(a_{n})_{n\geq1}$, $a_{n}\in T(\mathbb{Q})$, such that
$a_{n}=(a_{mn})^{m}$. In particular, $a_{n}$ is infinitely divisible. If
$T=(\mathbb{G}_{m})_{L/\mathbb{Q}}$, then $T(\mathbb{Q})=L^{\times}$, and
$\cap L^{\times m}=1$. Every torus $T$ can be embedded in a product of tori of
the form $(\mathbb{G}_{m})_{L/\mathbb{Q}}$, and so again $\cap T(\mathbb{Q}%
)^{m}=1$.} For any torus $T$ over $\mathbb{Q}{}$ such that $T(\mathbb{R}{})$
is compact, $\tilde{T}(\mathbb{\mathbb{A}{}}_{f})=0$ (ibid.
3.21).\footnote{First, $T(\mathbb{Q}{})$ is discrete in $T(\mathbb{A}{})$: let
$k$ be a finite splitting field for $T$; then $T(k)$ is discrete in
$T(\mathbb{A}_{k}{})$ by algebraic number theory; now use that $T(\mathbb{Q}%
{})=T(k)\cap T(\mathbb{A}{}_{\mathbb{Q}{}}).$
\par
The hypothesis implies that $T(\mathbb{A}{}_{f})/T(\mathbb{Q}{})$ is compact.
Therefore, for any choice of $\mathbb{Z}{}$-structure on $T$, the map
$T(\mathbb{\mathbb{\hat{Z}}{}})\rightarrow T(\mathbb{A}{}_{f})/T(\mathbb{Q}%
{})$ has finite kernel and cokernel. Now use that $\cap T(\mathbb{\hat{Z}%
)}^{m}=1$.}
\end{example}

\begin{remark}
\label{cc6}Let $T$ be a separable torus, $T=\varprojlim_{n}T_{n}$, over an
algebraically closed field $k$. Then the maps $T_{n+1}(k)\rightarrow T_{n}(k)$
are surjective, and so $\varprojlim^{1}T_{n}(k)=0$. Moreover, $H^{i}%
(k,T_{n})=0$ for all $n$ and all $i>0$, and so $H^{i}(k,T)=0$ for all $i>0$.
\end{remark}

\subsection{Ad\`{e}lic cohomology}

We now take $k=\mathbb{Q}{}$, so that $\Gamma=\Gal(\mathbb{Q}{}^{\text{al}%
}/\mathbb{Q}{})$. For a finite set $S$ of primes of $\mathbb{Q}{}$,
$\mathbb{A}{}^{S}$ is the restricted product of the $\mathbb{Q}{}_{l}$ for
$l\notin S$, and for a finite number field $L$, $\mathbb{A}{}_{L}^{S}%
=L\otimes_{\mathbb{Q}{}}\mathbb{A}^{S}$. When $S$ is empty, we omit it from
the notation.

For a torus $T$ over $\mathbb{Q}{}$, define
\[
H^{r}(\mathbb{A}{}^{S},T)=\dlim_{L}H^{r}(\Gamma_{L/\mathbb{Q}{}}%
,T({}\mathbb{A}{}_{L}^{S}))
\]
(limit over the finite Galois extensions $L$ of $\mathbb{Q}{}$ contained in
$\mathbb{Q}{}^{\text{al}}{}$).

\begin{proposition}
\label{ca3}Let $T$ be a torus over $\mathbb{Q}{}$. For all $r\geq1$, there is
a canonical isomorphism
\[
H^{r}(\mathbb{A}{}^{S},T)\cong\oplus_{l\notin S}H^{r}(\mathbb{Q}{}%
_{l},T)\text{\thinspace.}%
\]

\end{proposition}

\begin{proof}
Let $L/\mathbb{Q}{}$ be a finite Galois extension. For each prime $l$ of
$\mathbb{Q}{}$, choose a prime $v$ of $L$ lying over it, and set $L^{l}=L_{v}%
$. Then
\[
H^{r}(\Gamma_{L^{l}/\mathbb{Q}{}_{l}},T(L^{l}))\overset{\text{df}}{=}%
H^{r}(\Gamma_{L_{v}/\mathbb{Q}{}_{l}},T(L_{v}))
\]
is independent of the choice of $v$ up to a canonical isomorphism (i.e., it is
well-defined). Moreover,%
\[
H^{r}(\Gamma_{L/\mathbb{Q}{}},T(\mathbb{A}{}_{L}^{S}))\cong\oplus_{l\notin
S}H^{r}(L^{l}/\mathbb{Q}{}_{l},T(L^{l}))
\]
(\cite{platonov1994}, Proposition 6.7, p298). Now pass to the direct limit
over $L$.
\end{proof}

\begin{corollary}
\label{ca3m}A short exact sequence%
\[
0\rightarrow T^{\prime}\rightarrow T\rightarrow T^{\prime\prime}\rightarrow0
\]
of tori gives rise to a long exact sequence%
\[
\cdots\rightarrow H^{r}(\mathbb{A}^{S},T^{\prime})\rightarrow H^{r}%
(\mathbb{A}{}^{S},T)\rightarrow H^{r}(\mathbb{A}{}^{S},T^{\prime\prime
})\rightarrow H^{r+1}(\mathbb{A}{}^{S},T^{\prime})\rightarrow\cdots.
\]

\end{corollary}

\begin{proof}
Take the direct sum of the cohomology sequences over the $\mathbb{Q}{}_{l}$
and apply the proposition.
\end{proof}

For a torus $T$, let%
\[
T(\mathbb{\bar{A}}^{S}\mathbb{)=}{}\dlim_{L}T(\mathbb{A}{}_{L}^{S})\text{,}%
{}{}\quad\text{(limit over }L\subset\mathbb{Q}{}^{\text{al}}\text{ with
}[L\colon\mathbb{Q}{}]<\infty\text{),}%
\]
and, for a separable protorus $T$, let%
\[
T(\mathbb{\bar{A}}^{S}\mathbb{)}\mathbb{=}\varprojlim Q({}\mathbb{\bar{A}}%
^{S}\mathbb{)}\text{,}\quad\text{(limit over the algebraic quotients of
}T\text{)}.
\]
Endow each $Q(\mathbb{\bar{A}}^{S}\mathbb{)}{}$ with the discrete topology and
$T(\mathbb{\bar{A}}^{S}\mathbb{)}{}$ with the inverse limit topology, and
define%
\[
H^{r}(\mathbb{A}{}^{S},T)=H^{r}(\Gamma,T(\mathbb{\bar{A}}^{S})){}%
\]
where $H^{r}(\Gamma,T(\mathbb{\bar{A}}^{S}))$ is computed using continuous
cochains (profinite topology on $\Gamma$). When $T$ is a torus, this coincides
with the previous definition.

\begin{proposition}
\label{ca4}Let $T$ be a separable protorus over $\mathbb{Q}$.

\begin{enumerate}
\item There is a canonical homomorphism
\[
H^{1}(\mathbb{Q}{},T)\rightarrow H^{1}(\mathbb{A}{}^{S},T)\text{.}%
\]

\item Write $T$ as the limit of a strict inverse system of tori,
$T=\varprojlim T_{n}$. For each $r\geq0$, there is a canonical exact sequence
\[
0\rightarrow\plim{}^{1}T_{n}(\mathbb{A}{}^{S})\rightarrow H^{1}(\mathbb{A}%
{}^{S},T)\rightarrow\plim{}H^{1}(\mathbb{A}{}^{S},T_{n})\rightarrow0\text{.}%
\]

\end{enumerate}
\end{proposition}

\begin{proof}
(a) For each algebraic quotient $Q$ of $T$, $Q(\mathbb{Q}{}^{\text{al}%
})\hookrightarrow Q(\mathbb{\bar{A}}^{S}\mathbb{)}{}$ is continuous (both
groups are discrete), and hence the inverse limit $T(\mathbb{Q}{}^{\text{al}%
})\rightarrow T(\mathbb{\bar{A}}^{S}\mathbb{)}{}$ of these homomorphisms is continuous.

(b) The map $T_{n+1}(\mathbb{\bar{A}}^{S})\rightarrow T_{n}(\mathbb{\bar{A}}%
{}^{S}){}$ is surjective (in fact, $T_{n+1}(\mathbb{A}{}_{L}^{S})\rightarrow
T_{n}(\mathbb{A}{}_{L}^{S})$ will be surjective once $L$ is large enough to
split $T$), and so the proof of (\ref{cc2}) applies.
\end{proof}

\begin{remark}
\label{ca5}For any finite set $S$ of primes of $\mathbb{Q}{}$,
\[
T(\mathbb{\bar{A}}{})\cong T(\mathbb{\bar{A}}^{S})\times T({}(\prod_{l\in
S}\mathbb{Q}{}_{l})\otimes\mathbb{Q}{}^{\text{al}})
\]
as a topological group, and so
\[
H^{r}(\mathbb{A}{},T)\cong H^{r}(\mathbb{A}{}^{S},T)\times%
{\textstyle\prod\nolimits_{l\in S}}
H^{r}(\mathbb{Q}_{l},T)\text{.}%
\]

\end{remark}

\begin{remark}
\label{ca6}Let $T=\varprojlim(T_{n},u_{n})$ be a separable pro-torus over
$\mathbb{Q}{}$. It may happen that each $T_{n}$ satisfies the Hasse principle
but $T$ does not. In this case, we get a diagram
\[
\begin{CD}
0 @>>> \plim{}_n^{1}T_{n}(\mathbb{Q}) @>>> H^{1}(\mathbb{Q},T)
@>>> \plim_n H^{1}(\mathbb{Q},T_{n}) @>>> 0 \\
@.@VV{a}V @VV{b}V @VV{c}V @. \\
0 @>>> \plim{}_n^{1}T_{n}(\mathbb{A}) @>>> H^{1}(\mathbb{A},T) @>>>
\plim{}_n H^{1}(\mathbb{A},T_{n}) @>>> 0
\end{CD}
\]
in which $c$ is injective and
\[
\Ker(a)\cong\Ker(b)\neq0\text{.}%
\]
Let $(t_{n})$ be an element of $\prod T_{n}(\mathbb{Q}{})$ that is not in the
image of $1-u$ on $\prod T_{n}(\mathbb{Q}{})$ but is in the image of $1-u$ on
$\prod T_{n}(\mathbb{A}{})$. Then
\[
P=\plim{}(T_{n},u_{n}\cdot t_{n-1})
\]
is a nontrivial $T$-torsor under over $\mathbb{Q}{}$ that becomes trivial over
$\mathbb{A}{}$.\newpage
\end{remark}

\section{The cohomology of the Serre and Weil-number protori.}

Throughout this section, $\mathbb{Q}{}^{\text{al}}$ is the algebraic closure
of $\mathbb{Q}{}$ in $\mathbb{C}{}$, $\Gamma=\Gal(\mathbb{Q}{}^{\text{al}%
}/\mathbb{Q}{})$, and $\mathcal{F}{}$ is the set of CM-subfields $K$ of
$\mathbb{Q}{}^{\text{al}}$, finite and Galois over $\mathbb{Q}{}$.

\subsection{The Serre protorus $S$.}

For $K\in\mathcal{F}{}$, the Serre group $S^{K}$ for $K$ is the quotient of
$(\mathbb{G}_{m})_{K/\mathbb{Q}{}}$ such that%
\[
X^{\ast}(S^{K})=\left\{  f\colon\Gamma_{K/\mathbb{Q}{}}\rightarrow\mathbb{Z}%
{}\mid f(\sigma)+f(\iota\sigma)\text{ is constant}\right\}  .
\]
The constant value of $f(\sigma)+f(\iota\sigma)$ is called the weight of $f$.
There is an exact sequence
\begin{equation}
0\rightarrow(\mathbb{G}_{m})_{K_{+}/\mathbb{Q}{}}%
\xr{(\text{incl.},\Nm_{F/ \mathbb{Q}}^{-1})}(\mathbb{G}_{m})_{K/\mathbb{Q}{}%
}\times\mathbb{G}_{m}\rightarrow S^{K}\rightarrow0 \label{e15}%
\end{equation}
of tori over $\mathbb{Q}{}$ corresponding to the obvious exact sequence of
character groups. The maps $\Nm_{K^{\prime}/K}\times\id$ induce homomorphisms
$S^{K^{\prime}}\rightarrow S^{K}$, and the Serre group is defined to be%
\[
S=\plim_{\mathcal{F}{}}S^{K}.
\]

\begin{lemma}
\label{cg0}There is an exact sequence%
\[
0\rightarrow H^{1}(\mathbb{Q}{},S^{K})\rightarrow\Br(K_{+})\rightarrow
\Br(K)\oplus\Br(\mathbb{Q}{})\rightarrow H^{2}(\mathbb{Q}{},S^{K}%
)\rightarrow0\text{.}%
\]

\end{lemma}

\begin{proof}
Except for the zero at right, the statement follows from (\ref{e15}) and
Hilbert's Theorem 90, but a theorem of Tate (\cite{milne1986d}, I 4.21, p80)
shows that
\[
H^{3}(K_{+},\mathbb{G}_{m})\cong\oplus_{v\text{ real}}H^{3}((K_{+}%
)_{v},\mathbb{G}_{m}),
\]
and $H^{3}(\mathbb{R}{},\mathbb{G}_{m})\cong H^{1}(\mathbb{R}{},\mathbb{G}%
_{m})=0$ (periodicity of the cohomology of finite cyclic groups).
\end{proof}

\begin{proposition}
\label{cg1}For any $K$ as above${}$,
\begin{align*}
H^{1}(\mathbb{Q}{},S^{K})  &  \overset{\cong}{\rightarrow}H^{1}%
(\mathbb{\mathbb{A}{}},S^{K}),\\
H^{2}(\mathbb{Q}{},S^{K})  &  \hookrightarrow H^{2}(\mathbb{A},S^{K})\text{.}%
\end{align*}

\end{proposition}

\begin{proof}
Apply the snake lemma to the diagram,%
\[
\begin{CD}%
0 @>>> \Br(K_+) @>>> \oplus_l\Br((K_+)_l)@>>>
\mathbb{Q}/\mathbb{Z}@>>>0\\
@.@VVV@VVV@VV{x\mapsto(2x,x)}V@.\\
0 @>>> \Br(K)\oplus\Br(\mathbb{Q}) @>>> \oplus_l\Br(K_l)\oplus\oplus_l\Br(\mathbb{Q}_l)
@>>>
\mathbb{Q}/\mathbb{Z}\oplus\mathbb{Q}/\mathbb{Z}@>>>0
\end{CD}
\]
in which the subscript $l$ denotes $-\otimes\mathbb{Q}{}_{l}$. For the
vertical map at right we have used that
\begin{align*}
\inv_{E}\circ\Res  &  =[E\colon F]\cdot\inv_{F}\\
\inv_{F}\circ\mathrm{Cor}  &  =\inv_{E}%
\end{align*}
for an extension $F\subset E$; see \cite{serre1962}, XI \S 2, Proposition 1.
\end{proof}

\begin{proposition}
\label{cg2}There is a canonical exact sequence%
\[
0\rightarrow\plim{}_{\mathcal{F}{}}C(K)\rightarrow\plim{}_{\mathcal{F}{}}%
^{1}S^{K}(\mathbb{Q}{})\rightarrow\plim{}_{\mathcal{F}{}}^{1}S^{K}%
(\mathbb{A}{})\rightarrow0
\]
and a commutative diagram%
\[%
\begin{array}
[c]{ccccccc}%
0 & \rightarrow & \mathbb{Q}{}^{\times} & \rightarrow & S(\mathbb{Q}{}) &
\rightarrow & \plim{}^{1}K_{+}^{\times}\\
&  & \downarrow &  & \downarrow &  & \downarrow\\
0 & \rightarrow & \mathbb{A}{}^{\times}\times\varprojlim\frac{\left(
K_{+}\otimes\mathbb{C}{}\right)  ^{\times}}{\left(  K_{+}\otimes
\mathbb{\mathbb{R}{}}{}\right)  ^{\times}} & \rightarrow & S(\mathbb{A}{}) &
\rightarrow & \plim{}^{1}\mathbb{I}(K_{+})
\end{array}
\]

\end{proposition}

\begin{proof}
From (\ref{e15}), we obtain an inverse system of exact sequences
\[
0\rightarrow K_{+}^{\times}\rightarrow K^{\times}\times\mathbb{Q}{}^{\times
}\rightarrow S^{K}(\mathbb{Q}{})\rightarrow0,
\]
and hence an exact sequence (\ref{il03}, \ref{il04}),
\begin{equation}
0\rightarrow\mathbb{Q}{}^{\times}\rightarrow S(\mathbb{Q}{})\rightarrow
\plim{}^{1}K_{+}^{\times}\rightarrow\plim{}^{1}K^{\times}\rightarrow
\plim{}^{1}S^{K}(\mathbb{Q}{})\rightarrow0. \label{1e}%
\end{equation}
Similarly, there is an exact sequence
\begin{equation}
0\rightarrow\mathbb{A}{}^{\times}\times\varprojlim\frac{\left(  K_{+}%
\otimes\mathbb{C}{}\right)  ^{\times}}{\left(  K_{+}\otimes\mathbb{\mathbb{R}%
{}}{}\right)  ^{\times}}\rightarrow S(\mathbb{A}{})\rightarrow\plim{}^{1}%
\mathbb{I}{}(K_{+})\rightarrow\plim{}^{1}\mathbb{I}{}(K)\rightarrow
\plim{}^{1}S^{K}(\mathbb{A}{})\rightarrow0\text{.} \label{1f}%
\end{equation}
Thus, there are isomorphisms%
\begin{align*}
\plim{}^{1}\left(  K^{\times}/K_{+}^{\times}\right)   &  \cong\plim{}^{1}%
S^{K}(\mathbb{Q}{})\\
\plim{}^{1}\left(  \mathbb{I}{}(K)/\mathbb{I}{}(K_{+})\right)   &
\cong\plim{}^{1}S^{K}(\mathbb{A}{})\text{,}%
\end{align*}
and so the exact sequence follows from Proposition \ref{il10}. The commutative
diagram is obvious.
\end{proof}

\subsection{The Weil-number protorus $P$}

We now fix a $p$-adic prime $w$ of $\mathbb{Q}{}^{\text{al}}$, and we write
$w_{K}$ (or $w$) for the prime it induces on a subfield $K$ of $\mathbb{Q}%
^{\text{al}}$. The completion $(\mathbb{Q}{}^{\text{al}})_{w}$ of
$\mathbb{Q}{}^{\text{al}}$ at $w$ is algebraically closed, and we let
$\mathbb{\mathbb{Q}{}}_{p}^{\text{al}}$ denote the algebraic closure of
$\mathbb{Q}{}_{p}$ in $(\mathbb{Q}{}^{\text{al}})_{w}$. Its residue field,
which we denote $\mathbb{F}{}$, is an algebraic closure of $\mathbb{F}{}_{p}$.

A \emph{Weil }$p^{n}$\emph{-number} is an algebraic number $\pi$ for which
there exists an integer $m$ (the \emph{weight }of $\pi$) such that $\rho
\pi\cdot\overline{\rho\pi}=(p^{n})^{m}$ for all $\rho\colon\mathbb{Q}{}%
[\pi]\rightarrow\mathbb{C}{}$ . Let $W(p^{n})$ be the set of all Weil $p^{n}%
$-numbers in $\mathbb{Q}{}^{\text{al}}$. It is an abelian group, and for
$n|n^{\prime}$, $\pi\mapsto\pi^{n^{\prime}/n}$ is a homomorphism
$W(p^{n})\rightarrow W(p^{nn^{\prime}})$. Define%
\[
W=\dlim W(p^{n})\text{.}%
\]
There is an action of $\Gal(\mathbb{Q}{}^{\text{al}}/\mathbb{Q}{})$ on $W$,
and $P$ is defined to be the protorus over $\mathbb{Q}{}$ such that%
\[
X^{\ast}(P)=W\text{.}%
\]

For $\pi\in W(p^{n})$ and a $p$-adic prime $v$ of a finite number field
containing $\pi$, define%
\[
s_{\pi}(v)=\frac{\ord_{v}(\pi)}{\ord_{v}(p^{n})}\text{.}%
\]
Then $s_{\pi}(v)$ is well-defined for $\pi\in W$, i.e., it does not depend on
the choice of a representative of $\pi$.

Let $K$ be a CM field in $\mathbb{Q}{}^{\text{al}}$, finite and Galois over
$\mathbb{Q}{}$. Define $W^{K}$ to be the set of $\pi\in W$ having a
representative in $K$ and such that%
\[
f_{\pi}^{K}(v)\overset{\text{df}}{=}s_{\pi}(v)\cdot\lbrack K_{v}%
\colon\mathbb{Q}{}_{p}]\in\mathbb{Z}{}%
\]
for all $p$-adic primes of $K$, and define $P^{K}$ to be the torus over
$\mathbb{Q}{}$ such that%
\[
X^{\ast}(P^{K})=W^{K}.
\]
Then%
\[
W=\dlim_{\mathcal{F}{}}W^{K},\quad P=\varprojlim_{\mathcal{F}{}}P^{K}.
\]

Let $X$ and $Y$ respectively be the sets of $p$-adic primes of $K$ and $K_{+}%
$. Then (e.g., \cite{milne2001a}, A.6), there is an exact sequence%
\begin{equation}
\begin{CD} 0@>>>W^K@>{\pi\mapsto}>{(f^K_{\pi},wt(\pi))}>\mathbb{Z}^X\times\mathbb{Z} @>{(f,m)\mapsto}>{f|Y-n(w_K)\cdot m\cdot \sum_{v\in Y}v}>\mathbb{Z}^Y@>>>0\\ \end{CD} \label{e18}%
\end{equation}
where $n(w_{K})$ is the local degree $[K_{w}:\mathbb{Q}{}_{p}]$. Our fixed
$p$-adic prime $w$ allows us to identify $X$ with $\Gamma_{K/\mathbb{Q}{}%
}/D(w_{K})$ and $Y$ with $\Gamma_{K_{+}/\mathbb{Q}{}}/D(w_{K_{+}})$.
Therefore, (\ref{e18}) is the sequence of character groups of an exact
sequence%
\begin{equation}
\begin{CD} 0@>>>(\mathbb{G}_{m})_{K(w)_{+}/\mathbb{Q}} @>\left(\begin{smallmatrix} \text{incl.}\\ \Nm^{-n(w_K)}\end{smallmatrix}\right)>> (\mathbb{G}_{m})_{K(w)/\mathbb{Q}{}}\times\mathbb{G}_{m} @>>> P^{K}@>>>0. \end{CD} \label{e16}%
\end{equation}

Let $K^{\prime}\supset K$ with $K^{\prime}\in\mathcal{F}{}$. Then
$W^{K}\subset W^{K^{\prime}}$, and there is map $X^{\prime}\rightarrow X$ from
the $p$-adic primes of $K^{\prime}$ to those of $K$. If $v^{\prime}\mapsto v$,
then $s_{\pi}(v^{\prime})=s_{\pi}(v)$, and so $f_{\pi}(v^{\prime}%
)=[K_{w}^{\prime}\colon K_{w}]\cdot f_{\pi}(v)$. Therefore, the diagram%
\[
\begin{CD}
0 @>>> W^{K} @>>> \mathbb{Z}^{X}\times\mathbb{Z} @>>> \mathbb{Z}^{Y} @>>>0\\
@.\cap@.@VV{a\times\id}V@VVV\\
0 @>>> W^{K^{\prime}} @>>> \mathbb{Z}{}^{X^{\prime}}\times\mathbb{Z}@>>>\mathbb{Z}^{Y^{\prime}}@>>> 0
\end{CD}
\]
commutes with $a$ equal to $[K_{w}^{\prime}\colon K_{w}]\times$ the map
induced by $Y^{\prime}\rightarrow Y$. Therefore,
\begin{equation}
\begin{CD}0 @>>> (\mathbb{G}_{m})_{K^{\prime}(w)_{+}/\mathbb{Q}} @>>> (\mathbb{G}_{m})_{K^{\prime}(w)/\mathbb{Q}}\times\mathbb{G}_{m} @>>> P^{K^{\prime}} @>>> 0\\ @. @VV{b}V@VV{\Nm_{K^{\prime}/K}^{[K_{w}^{\prime}\colon K_{w}]}\times\id}V@VV{c}V@.\\ 0 @>>> (\mathbb{G}_{m})_{K(w)_{+}/\mathbb{Q}{}} @>>> (\mathbb{G}_{m})_{K(w)/\mathbb{Q}}\times\mathbb{G}_{m} @>>> P^{K} @>>> 0. \end{CD} \label{e25}%
\end{equation}
commutes.

\begin{lemma}
\label{cg3}There is an exact sequence%
\[
0\rightarrow H^{1}(\mathbb{Q}{},P^{K})\rightarrow\Br(K(w)_{+})\rightarrow
\Br(K(w))\oplus\Br(\mathbb{Q}{})\rightarrow H^{2}(\mathbb{Q}{},P^{K}%
)\rightarrow0.
\]

\end{lemma}

\begin{proof}
Same as that of Lemma \ref{cg0}.
\end{proof}

\begin{proposition}
\label{cg4}Let $K\in\mathcal{F}{}$.

\begin{enumerate}
\item If the local degree $n(w_{K})=[K_{w}\colon\mathbb{Q}{}_{p}]$ is even,
then there is an exact sequence%
\[
0\rightarrow H^{1}(\mathbb{Q}{},P^{K})\rightarrow H^{1}(\mathbb{A}{}%
,P^{K})\rightarrow\frac{1}{2}\mathbb{Z}{}/\mathbb{Z}{}\rightarrow0;
\]
otherwise,%
\[
H^{1}(\mathbb{Q}{},P^{K})\cong H^{1}(\mathbb{A},P^{K})
\]

\item The map%
\[
H^{2}(\mathbb{Q}{},P^{K})\rightarrow H^{2}(\mathbb{A},P^{K})
\]
is injective.
\end{enumerate}
\end{proposition}

\begin{proof}
If complex conjugation $\iota\in D(w_{K})$, then $P^{K}=\mathbb{G}_{m}$ and
the statement is obvious. Thus, we may assume $\iota\notin D(w_{K})$. Consider
the diagram%
\[
\begin{CD}%
0 @>>> \Br(K(w)_+) @>>> \oplus_l\Br((K(w)_+)_l)@>>>
\mathbb{Q}/\mathbb{Z}@>>>0\\
@.@VVV@VVV@VV{x\mapsto(2x,n(w_K)x)}V@.\\
0 @>>> \Br(K(w))\oplus\Br(\mathbb{Q}) @>>> \oplus_l\Br(K(w)_{l})\oplus
\oplus_l\Br(\mathbb{Q}_l)@>>>
\mathbb{Q}/\mathbb{Z}\oplus\mathbb{Q}/\mathbb{Z}@>>>0.
\end{CD}
\]
When $n(w_{K})$ is odd, $(2,n(w_{K}))$ is injective and the snake lemma shows
that that $H^{1}(\mathbb{Q}{},P^{K})\cong H^{1}(\mathbb{A},P^{K})$. When
$n(w_{K})$ is even, the sequence of kernels is%
\[
0\rightarrow\Br(K(w)/K(w)_{+})\rightarrow\oplus_{l}\Br(K(w)_{l}/(K(w)_{+}%
)_{l})\rightarrow\frac{1}{2}\mathbb{Z}{}/\mathbb{Z}{}\rightarrow0,
\]
which class field theory shows to be exact. Again, the snake lemma gives the result.
\end{proof}

\begin{proposition}
\label{cg5}The canonical map
\[
\plim{}_{\mathcal{F}{}}^{1}P^{K}(\mathbb{Q}{})\rightarrow\plim{}_{\mathcal{F}%
{}}^{1}P^{K}(\mathbb{A}{})
\]
is an isomorphsim.
\end{proposition}

\begin{proof}
As in the proof of \ref{cg2}, there are canonical isomorphism
\begin{align*}
\plim{}^{1}\left(  K(w)^{\times}/K(w)_{+}^{\times}\right)   &  \cong
\plim{}^{1}P^{K}(\mathbb{Q}{})\\
\plim{}^{1}\left(  \mathbb{I(}{}K(w)/\mathbb{I}({}K(w)_{+})\right)   &
\cong\plim{}^{1}P^{K}(\mathbb{A})
\end{align*}
and so the statement follows from Proposition \ref{il12} (\ref{e24}).
\end{proof}

\begin{lemma}
\label{cg5n}For each $K\in\mathcal{F}{}$, there exists an $K^{\prime}%
\in\mathcal{F}{}$ containing $K$ for which the map%
\[
H^{1}(\mathbb{Q}{},P^{K^{\prime}})\rightarrow H^{1}(\mathbb{Q}{},P^{K})
\]
is zero.
\end{lemma}

\begin{proof}
After possibly enlarging $K$, we may assume that $\iota\notin D(w_{K})$. The
map $b$ in (\ref{e25}) is $\Nm_{K^{\prime}(w)_{+}/K(w)_{+}}^{[K_{+w}^{\prime
}\colon K_{+w}]}$. \noindent If $K^{\prime}$ is chosen so that $2|[K^{\prime
}(w)\colon K(w)]$, then $H^{2}(b)$ is zero on the kernel of $\Br(K^{\prime
}(w)_{+})\rightarrow\Br(K(w))$, and hence $H^{1}(c)$ is zero.
\end{proof}

\begin{lemma}
\label{cg5m}The groups $\varprojlim_{\mathcal{F}{}}H^{1}(\mathbb{Q}{},P^{K})$
and $\varprojlim_{\mathcal{F}{}}^{1}{}H^{1}(\mathbb{Q}{},P^{K})$ are both zero.
\end{lemma}

\begin{proof}
Lemma \ref{cg5n} shows that $(H^{1}(\mathbb{Q}{},P^{K}))_{K\in\mathcal{F}{}}$
admits a cofinal subsystem in which the transition maps are zero, and so the
map $u$ in (\ref{e10}) is zero.
\end{proof}

\begin{lemma}
\label{cg5s}The groups $\varprojlim_{\mathcal{F}{}}H^{1}(\mathbb{A},P^{K})$
and $\varprojlim_{\mathcal{F}{}}^{1}{}H^{1}(\mathbb{A},P^{K})$ are both zero.
\end{lemma}

\begin{proof}
For $K$ sufficiently large, (\ref{cg4}) shows that $H^{1}(\mathbb{Q}{}%
,P^{K})\cong H^{1}(\mathbb{A}{},P^{K})$ and so this follows from the previous lemma.
\end{proof}

\begin{proposition}
\label{cg5p}There are canonical isomorphisms%
\begin{align}
\plim{}^{1}P^{K}(\mathbb{Q}{})  &  \cong H^{1}(\mathbb{Q}{},P)\label{e34}\\
H^{2}(\mathbb{Q}{},P)  &  \cong\varprojlim H^{2}(\mathbb{Q}{},P^{K}%
)\label{e35}\\
\plim{}^{1}P^{K}(\mathbb{A})  &  \cong H^{1}(\mathbb{A}{},P)\label{e36}\\
H^{2}(\mathbb{A},P)  &  \cong\varprojlim H^{2}(\mathbb{A},P^{K})\text{.}
\label{e37}%
\end{align}

\end{proposition}

\begin{proof}
Combine Proposition \ref{cc2} with Lemmas \ref{cg5m} and \ref{cg5s}$.$
\end{proof}

\begin{remark}
\label{cg5q}Assume $\iota\notin D(w_{K})$. The argument in the proof of Lemma
\ref{cg5n} shows that, if the local degree $[L_{w}\colon K_{w}]$ is even, then
all the vertical maps in the diagram%
\[
\begin{CD}
0 @>>> H^{1}(\mathbb{Q},P^{L}) @>>> H^{1}(\mathbb{A},P^{L}) @>>> \frac{1}{2}\mathbb{Z}/\mathbb{Z} @>>> 0\\
@.@VVV@VVV@VVV@.\\
0 @>>> H^{1}(\mathbb{Q}{},P^{K}) @>>> H^{1}(\mathbb{A}%
{},P^{K}) @>>> \frac{1}{2}\mathbb{Z}{}/\mathbb{Z} @>>> 0
\end{CD}
\]
are zero.
\end{remark}

\begin{proposition}
\label{cg5r}There are canonical isomorphisms%
\begin{align*}
P(\mathbb{Q}{})  &  \cong\mathbb{Q}{}^{\times}\\
P(\mathbb{A}{}_{f})  &  \cong\mathbb{A}{}_{f}^{\times}.
\end{align*}

\end{proposition}

\begin{proof}
For $K$ sufficiently large, $P^{K}\cong P_{0}^{K}\oplus\mathbb{G}_{m}$ where
$X^{\ast}(P_{0}^{K})$ consists of the Weil numbers of weight $0$. Thus $P\cong
P_{0}\oplus\mathbb{G}_{m}$. For each sufficiently large $K$, the projection
$P_{0}\rightarrow P_{0}^{K}$ factors into $P_{0}\twoheadrightarrow
\widetilde{P_{0}^{K}}\rightarrow P_{0}^{K}$ where $\widetilde{P_{0}^{K}}$ is
the universal covering of $P_{0}^{K}$. Therefore, $P_{0}=\varprojlim
_{K}\widetilde{P_{0}^{K}}$, and $P_{0}(\mathbb{Q}{})=\varprojlim
\widetilde{P_{0}^{K}}(\mathbb{Q}{})$, which is $0$ by (\ref{cc5}). Similarly,
$P_{0}(\mathbb{A}_{f}{})=\varprojlim\widetilde{P_{0}^{K}}(\mathbb{A}_{f})=0$.
\end{proof}

\subsection{The cohomology of $S/P$.}

Let $K\in\mathcal{F}{}$, and assume $\iota\notin D(w_{K})$. Let $\mathfrak{p}%
{}$ be the prime ideal in $\mathcal{O}{}_{K}$ corresponding to $w_{K}$. For
some $h$, $\mathfrak{p}{}^{h}$ will be principal, say $\mathfrak{p}{}^{h}%
=(a)$. Let $\alpha=a^{2n}$ where $n=(U(K)\colon U(K_{+}))$. Then, for $f\in
X^{\ast}(S^{K})$, $f(\alpha)$ is independent of the choice of $a$, and it is a
Weil $p^{2nf(\mathfrak{p}{}/p)}$-number of weight $wt(f)$. The map $f\mapsto
f(\alpha)\colon X^{\ast}(S^{K})\rightarrow W^{K}$ is a surjective homomorphism
(e.g., \cite{milne2001a}, A.8). Thus, it corresponds to an injective
homomorphism $\rho^{K}\colon P^{K}\rightarrow S^{K}$ which can also be
characterized as the unique homomorphism rendering%
\[
\begin{CD}
0 @>>> (\mathbb{G}_{m})_{K(w)_{+}/\mathbb{Q}{}} @>>> (\mathbb{G}%
_{m})_{K(w)/\mathbb{Q}{}}\times\mathbb{G}_{m} @>>> P^{K} @>>> 0\\
@.@VV{\text{incl.}}V@VV{\text{incl.}\times\id}V@VV{\rho_K}V@.\\
0 @>>> (\mathbb{G}_{m})_{K_{+}/\mathbb{Q}} @>>>
(\mathbb{G}_{m})_{K/\mathbb{Q}{}}\times\mathbb{G}_{m} @>>> S^{K} @>>> 0
\end{CD}
\]
commutative. For varying $K$, the $\rho^{K}$ define a morphism of inverse
systems. Therefore, on passing to the inverse limit, we obtain an injective
homomorphism $\rho\colon P\rightarrow S$ of protori.

\begin{proposition}
\label{cg6}The map
\[
H^{1}(\mathbb{Q},S^{K}/P^{K})\rightarrow H^{1}(\mathbb{A}{},S^{K}/P^{K})
\]
is injective on the image of $H^{1}(\mathbb{Q},S^{L}/P^{L})$ for any $L\supset
K$ such that the local degree $[L_{w}\colon K_{w}]$ at $w$ is even. Therefore,
the map
\[
\varprojlim H^{1}(\mathbb{Q},S^{K}/P^{K})\rightarrow\varprojlim H^{1}%
(\mathbb{A}{},S^{K}/P^{K})
\]
is injective.
\end{proposition}

\begin{proof}
Diagram chase in
\[
\begin{CD}
H^1(\mathbb{Q},P^K)@>>>H^{1}(\mathbb{Q},S^K) @>>> H^{1}(\mathbb{Q},S^K/P^K)
@>>> H^{2}(\mathbb{Q},P^K)\\
@VV(\ref{cg4})V@V{\cong}V{(\ref{cg1})}V@VVV@V\textrm{inj.}V(\ref{cg4})
V\\
H^1(\mathbb{A},P^K)@>>>H^{1}(\mathbb{A},S^K) @>>>H^{1}(\mathbb{A},S^K/P^K)
@>>>H^{2}(\mathbb{A},P^K),
\end{CD}
\]
using (\ref{cg5q}).
\end{proof}

\begin{proposition}
An element $a\in H^{1}(\mathbb{A},S^{K}/P^{K})$ arises from an element of
$H^{1}(\mathbb{Q}{},S^{K}/P^{K})$ if and only if its image in $H^{2}%
(\mathbb{A}{},P^{K})$ under the connecting homomorphism arises from an element
of $H^{2}(\mathbb{Q}{},P^{K})$.
\end{proposition}

\begin{proof}
Diagram chase in%
\[
\begin{CD}
H^{1}(\mathbb{Q},S^K) @>>> H^{1}(\mathbb{Q},S^K/P^K) @>>> H^{2}(\mathbb{Q},P^K) @>>> H^{2}(\mathbb{Q},S^K) \\
@V{\cong}V{(\ref{cg1})}V@VVV@VV
V@V\textrm{inj.}V(\ref{cg1})V\\
H^{1}(\mathbb{A},S^K) @>>>H^{1}(\mathbb{A},S^K/P^K)
@>>>H^{2}(\mathbb{A},P^K) @>>> H^{2}(\mathbb{A},S^K).
\end{CD}
\]

\end{proof}

\begin{proposition}
\label{cg7}There is a canonical exact sequence
\[
0\rightarrow\varprojlim_{\mathcal{F}{}}C(K)\rightarrow H^{1}(\mathbb{Q}%
{},S/P)\overset{a}{\rightarrow}H^{1}(\mathbb{A}{},S/P)\text{.}%
\]

\end{proposition}

\begin{proof}
Because $\varprojlim H^{1}(\mathbb{Q}{},S^{K}/P^{K})\rightarrow\varprojlim
H^{1}(\mathbb{A}{},S^{K}/P^{K})$ is injective,%
\[
\Ker(a)\cong\Ker(\plim{}^{1}(S^{K}/P^{K})(\mathbb{Q}{})\overset{b}%
{\rightarrow}\plim{}^{1}(S^{K}/P^{K})(\mathbb{A}));
\]
cf. (\ref{ca6}). On comparing (\ref{e15}) and (\ref{e16}), we obtain an exact
sequence%
\[
0\rightarrow(\mathbb{G}_{m})_{K(w)/\mathbb{Q}{}}/(\mathbb{G}_{m}%
)_{K(w)_{+}/\mathbb{Q}{}}\rightarrow(\mathbb{G}_{m})_{K/\mathbb{Q}{}%
}/(\mathbb{G}_{m})_{K_{+}/\mathbb{Q}{}}\rightarrow S^{K}/P^{K}\rightarrow
0\text{,}%
\]
which gives rise to an exact commutative diagram
\[
\begin{CD}%
\varprojlim^{1}(K(w)^{\times}/K(w)_{+}^{\times}) @>c>>
\varprojlim^{1}(K^{\times}/K_{+}^{\times}) @>>> \varprojlim
^{1}(S^{K}/P^{K})(\mathbb{\mathbb{Q}{}}) @>>> 0\\
@V\mathrm{surj.}V(\ref{il01m})V@VVdV@VVbV\\
\varprojlim^{1}(\mathbb{I}(K(w))/\mathbb{I}(K(w)_{+}^{\times}) @>>>
\varprojlim^{1}(\mathbb{I}{}(K)/\mathbb{I}{}(K_{+})) @>>>
\varprojlim^{1}(S^{K}/P^{K})(\mathbb{A}) @>>> 0.
\end{CD}
\]
The left-hand vertical map is surjective by \ref{il01m}, and the diagram gives
an exact sequence%
\[
\Ker(d\circ c)\overset{c}{\rightarrow}\Ker(d)\rightarrow\Ker(b)\rightarrow
0\text{.}%
\]
We can now apply Proposition \ref{il13}.
\end{proof}

\paragraph{Notes.}

Most of the calculations concerning the cohomology of $S^{K}$ and $P^{K}$ (but
not of $S$ or $P$ themselves) can be found already in \cite{langlands1987}.

\newpage

\section{The fundamental classes}

\subsection{The ad\`{e}lic fundamental class}

We use the Betti fibre functor $\omega_{B}$ to identify $\CM(\mathbb{Q}%
^{\text{al}})$ with $\Rep_{\mathbb{Q}{}}(S)$. Let $w_{\text{can}}$ denote the
cocharacter of $S$ such that $w_{\text{can}}(a)$ acts on an object of
$\CM(\mathbb{Q}^{\text{al}})$ of weight $m$ as $a^{m}$, and let $\mu
_{\text{can}}$ denote the cocharacter of $S_{\mathbb{Q}{}^{\text{al}}}$ such
that, for $f\in X^{\ast}(S^{K})$, $\langle\chi,f\rangle=f(i)$ where $i$ is the
given inclusion $K\hookrightarrow\mathbb{Q}{}^{\text{al}}$.

\subsubsection{The local fundamental class at $\infty$}

Let $\mathsf{R}_{\infty}$ (realization category at $\infty$) be the category
of pairs $(V,F)$ consisting of a $\mathbb{Z}{}$-graded finite-dimensional
complex vector space $V=\oplus_{m\in\mathbb{Z}{}}V^{m}$ and a semilinear
endomorphism $F$ such that $F^{2}=(-1)^{m}$. With the obvious tensor
structure, $\mathsf{R}_{\infty}$ becomes a Tannakian category with fundamental
group $\mathbb{G}_{m}$.

Let $(V,r)$ be a real representation of the Serre group $S$. Then
$w(r)=_{\text{df}}r\circ w_{\text{can}}$ defines a $\mathbb{Z}{}$-gradation on
$V\otimes\mathbb{C}{}$. Let $F$ be the map%
\[
v\mapsto r(\mu_{\text{can}}(i))\bar{v}\colon V\otimes\mathbb{C}{}\rightarrow
V\otimes\mathbb{C}{}\text{.}%
\]
Then $(V\otimes\mathbb{C}{},F)$ is an object of $\mathsf{R}_{\infty}$, and
$(V,r)\mapsto(V\otimes\mathbb{C}{},F)$ defines a tensor functor%
\[
\xi_{\infty}\colon\Rep_{\mathbb{R}}(S)\rightarrow\mathsf{R}_{\infty}.
\]
The functor $\xi_{\infty}$ defines a homomorphism $x_{\infty}\colon
\mathbb{G}_{m}\rightarrow S_{\mathbb{R}{}}$, which is equal to $w_{\text{can}%
/\mathbb{R}{}}$.

Let $\mathsf{R}_{\infty}^{\mathbb{G}_{m}}$ be the full subcategory of
$\mathsf{R}_{\infty}$ of objects of weight zero. For any $(V,F)$ of weight
zero,
\[
V^{F}=\{v\in V\mid Fv=v\}
\]
is a real form of $V$, and $(V,F)\mapsto V^{F}$ is a fibre functor on
$\mathsf{R}_{\infty}^{\mathbb{G}_{m}}$. Its composite with $\Rep_{\mathbb{R}%
{}}\left(  S\right)  ^{\mathbb{G}_{m}}\rightarrow\mathsf{R}_{\infty
}^{\mathbb{G}_{m}}$ is a fibre functor $\omega_{\infty}$, and $\wp
(\omega_{\infty})=_{\text{df}}\underline{\Hom}^{\otimes}(\omega_{B}%
,\omega_{\infty})$ is an $S/\mathbb{G}_{m}$-torsor. We define $c_{\infty}$ to
be the cohomology class of $\wp(\omega_{\infty})\wedge^{S/\mathbb{G}_{m}}S/P$.

\subsubsection{The local fundamental class at $p$.}

Let $\mathsf{R}_{p}$ (realization category at $p$) be the category of
$F$-isocrystals over $\mathbb{F}{}$. Thus, an object of $\mathsf{R}_{p}$ is a
pair $(V,F)$ consisting of a finite-dimensional vector space $V$ over
$B(\mathbb{F}{})$ and a semilinear isomorphism $F\colon V\rightarrow V$. With
the obvious tensor structure, $\mathsf{R}_{p}$ becomes a Tannakian category
with fundamental group $\mathbb{G}=\widetilde{\mathbb{G}_{m}}$.

There is a tensor functor%
\[
\xi_{p}\colon\CM(\mathbb{Q}^{\text{al}})\rightarrow\mathsf{R}_{p}%
\]
such that%
\[
\xi_{p}(h(A,e,m))=e\cdot H_{\text{crys}}^{\ast}(A_{\mathbb{F}{}})(m)
\]
where $A_{\mathbb{F}{}}$ is the reduction of $A$ at $w$. The functor $\xi_{p}$
defines a homomorphism $x_{p}\colon\mathbb{G}{}\rightarrow S_{/\mathbb{Q}%
{}_{p}}$ whose action on $X^{\ast}(S^{K})$ is $f\mapsto%
{\textstyle\sum_{\sigma\in D(w_{K})}}
f(\sigma)/(D(w_{K})\colon1)$.

Let $\mathsf{R}_{p}^{\mathbb{G}{}}$ be the full subcategory of $\mathsf{R}%
_{p}$ of objects of slope $0$. For any $(V,F)$ of slope zero,%
\[
V^{F}=\{v\in V\mid Fv=v\}
\]
is a $\mathbb{Q}{}_{p}$-form of $V$, and $(V,F)\mapsto V^{F}$ is a fibre
functor on $\mathsf{R}_{p}^{\mathbb{G}{}}$. Its composite with $\CM(\mathbb{Q}%
^{\text{al}})^{\mathbb{G}{}}\overset{\xi_{\infty}}{\rightarrow}\mathsf{R}%
_{p}^{\mathbb{G}}$ is a fibre functor on $\CM(\mathbb{Q}^{\text{al}%
})^{\mathbb{G}{}}$, and $\wp(\omega_{p})=_{\text{df}}\underline{\Hom}%
^{\otimes}(\omega_{B},\omega_{p})$ is a $S/\mathbb{G}{}$-torsor. We define
$c_{p}$ to be the cohomology class of $\wp(\omega_{p})\wedge^{S/\mathbb{G}{}%
}S/P$.

\subsubsection{The ad\`{e}lic fundamental class}

We define $c_{\mathbb{A}{}}\in H^{1}(\mathbb{A}{},S/P)$ to be the class
corresponding to $(0,c_{p},c_{\infty})$ under the isomorphism (\ref{ca5})
\[
H^{1}(\mathbb{A}{},S/P)\cong H^{1}(\mathbb{A}{}^{\{p,\infty\}},S/P)\times
H^{1}(\mathbb{Q}_{p},S/P)\times H^{1}(\mathbb{\mathbb{R}{}},S/P)\text{.}%
\]

\subsection{The global fundamental class}

\begin{lemma}
\label{fc06}Let $d_{p}^{K}$ and $d_{\infty}^{K}$ be the images of the classes
of $\mathsf{R}_{p}^{K}$ and $\mathsf{R}_{\infty}^{K}$ in $H^{2}(\mathbb{Q}%
{}_{p},P^{K})$ and $H^{2}(\mathbb{R},P^{K})$. There exists a unique element of
$H^{2}(\mathbb{Q}{},P^{K})$ with image%
\[
(0,d_{p}^{K},d_{\infty}^{K})\in H^{2}(\mathbb{A}^{\{p,\infty\}},P^{K})\times
H^{2}(\mathbb{Q}{}_{p},P^{K})\times H^{2}(\mathbb{R}{},P^{K})\text{.}%
\]

\end{lemma}

\begin{proof}
The uniqueness follows from (\ref{cg4}b). The existence is proved in
\cite{langlands1987} (also, \cite{milne1994a}, proof of 3.31).
\end{proof}

\begin{proposition}
\label{fc07}Let $d_{p}$ and $d_{\infty}$ be the images of the classes of
$\mathsf{R}_{p}$ and $\mathsf{R}_{\infty}$ in $H^{2}(\mathbb{Q}{}_{p},P)$ and
$H^{2}(\mathbb{R},P)$. There exists a unique element of $H^{2}(\mathbb{Q}%
{},P)$ with image%
\[
(0,d_{p},d_{\infty})\in H^{2}(\mathbb{A}^{\{p,\infty\}},P)\times
H^{2}(\mathbb{Q}{}_{p},P)\times H^{2}(\mathbb{R}{},P)\text{.}%
\]

\end{proposition}

\begin{proof}
Follows from the lemma, using (\ref{e35}) and (\ref{e37}).
\end{proof}

\begin{theorem}
\label{fc02}There exists a $c\in H^{1}(\mathbb{Q}{},S/P)$ mapping to
$c_{\mathbb{A}{}}\in H^{1}(\mathbb{A},S/P)$; any two such $c$'s have the same
image in $H^{1}(\mathbb{Q}{},S^{K}/P^{K})$ for all $K$; the set of such $c$'s
is a principal homogeneous space for $\varprojlim_{\mathcal{F}{}}C(K)$ where
$\mathcal{F}{}$ is the set of CM-subfields of $\mathbb{Q}{}^{\text{al}}$
finite over $\mathbb{Q}{}$.
\end{theorem}

We first prove three lemmas.

\begin{lemma}
\label{fc02m}Let $\mathsf{C}$ and $\mathsf{Q}$ be Tannakian categories with
commutative fundamental groups $G$ and $H$. Assume $\mathsf{C}$ is neutral.
Let $\xi\colon\mathsf{C}\rightarrow\mathsf{Q}$ be a tensor functor inducing an
injective homomorphism $H\rightarrow G$. Let $\omega_{\xi}$ denote the fibre
functor%
\[
\mathsf{C}^{H}\overset{\xi}{\rightarrow}\mathsf{Q}^{H}\xr{\Hom(\1,-)}\Vc_{k}%
\text{.}%
\]
For any fibre functor $\omega$ on $\mathsf{C}$, the class of the torsor
$\underline{\Hom}^{\otimes}(\omega,\omega_{\xi})$ in $H^{1}(k,G/H)$ maps to
the class of $\mathsf{Q}$ in $H^{2}(k,H)$ under the connecting homomorphism.
\end{lemma}

\begin{proof}
For a Tannakian category $\mathsf{T}$, write $\mathsf{T}^{\vee}$ for the gerbe
of fibre functors on $\mathsf{T}$. We are given a morphism $\xi^{\vee}%
\colon\mathsf{Q}^{\vee}\rightarrow\mathsf{C}^{\vee}$ bound by the injective
homomorphism $H\rightarrow G$ of commutative affine group schemes. Clearly
$\mathsf{Q}^{\vee}$ is the gerbe of local liftings of $\underline
{\Hom}^{\otimes}(\omega,\omega_{\xi})$ (\cite{giraud1971}, IV 2.5.4.1, p238),
and so its class is the image of $\underline{\Hom}^{\otimes}(\omega
,\omega_{\xi})$ under Giraud's definition of the connecting homomorphism
(ibid. IV 4.2), which coincides with the usual connecting homomorphism in the
commutative case (ibid. IV 3.4). Finally, the class of $\mathsf{Q}$ in
$H^{1}(k,H)$ is defined to be that represented by the gerbe $\mathsf{Q}^{\vee
}$.
\end{proof}

\begin{lemma}
\label{fc03}The image of $c_{\mathbb{A}{}}^{K}$ in $H^{2}(\mathbb{A}{},P^{K})$
under the connecting homomorphism arises from an element of $H^{2}%
(\mathbb{Q},P^{K})$.
\end{lemma}

\begin{proof}
Let $d_{p}^{K}$ and $d_{\infty}^{K}$ be the images of the classes of
$\mathsf{R}_{p}^{K}$ and $\mathsf{R}_{\infty}^{K}$ in $H^{2}(\mathbb{Q}{}%
_{p},P^{K})$ and $H^{2}(\mathbb{R},P^{K})$. According to Lemma \ref{fc02m}, we
have to prove that the element%
\[
(0,d_{p}^{K},d_{\infty}^{K})\in H^{2}(\mathbb{A}^{\{p,\infty\}},P^{K})\times
H^{2}(\mathbb{Q}{}_{p},P^{K})\times H^{2}(\mathbb{R}{},P^{K})
\]
arises from an element of $H^{2}(\mathbb{Q}{},P^{K})$, but this was shown in
Lemma \ref{fc06}.
\end{proof}

\begin{lemma}
\label{fc04}For every CM-subfield $K$ of $\mathbb{C}{}$, finite and Galois
over $\mathbb{\ Q}{}$, there exists a unique $c^{K}\in H^{1}(\mathbb{Q}%
{},S^{K}/P^{K})$ mapping to $c_{\mathbb{A}{}}^{K}\in H^{1}(\mathbb{A}%
,S^{K}/P^{K})$ and lifting to some $L\supset K$ with $[L_{w}\colon K_{w}]$ even.
\end{lemma}

\begin{proof}
Apply \ref{fc02}, \ref{fc03}.
\end{proof}

We now prove the theorem. Consider the diagram
\[
\begin{CD}
0 @>>> \plim{}^{1}(S^{K}/P^{K})(\mathbb{Q})
@>>> H^{1}(\mathbb{Q},S/P) @>>> \plim H^{1}(\mathbb{Q},S^{K}/P^{K}) \\
@.@VVbV @VVaV @VVV\\ 0 @>>> \plim{}^{1}(S^{K}/P^{K})(\mathbb{A}) @>>> H^{1}(\mathbb{A},S/P) @>>> \plim H^{1}(\mathbb{A},S^{K}/P^{K}),\\ @. @VVV @.@.\\ @. 0 \end{CD}
\]
The groups $S^{K}/P^{K}$ are anisotropic, which explains the $0$ at lower left
(\ref{il01m}). We have to show that there is a $c\in H^{1}(\mathbb{Q}{},S/P)$
mapping to $c_{\mathbb{A}{}}\in H^{1}(\mathbb{A}{},S/P)$. We know (\ref{fc04})
that, for each $K$, there is a unique element $c^{K}\in H^{1}(\mathbb{Q}%
{},S^{K}/P^{K})$ mapping to the image $c_{\mathbb{A}{}}^{K}$ of $c_{\mathbb{A}%
{}}$ in $H^{1}(\mathbb{A}{},S^{K}/P^{K})$ and lifting to some $L\supset K$
with $[L_{w}\colon K_{w}]$ even. Because of the uniqueness, the $c^{K}$ define
an element $\underleftarrow{c}\in$ $\plim H^{1}(\mathbb{Q}{},S^{K}/P^{K})$.
Choose $c\in H^{1}(\mathbb{Q}{},S/P)$ to map to $(c^{K})$. A diagram chase
using the surjectivity of $b$ shows that $c$ can be chosen to map to
$c_{\mathbb{A}{}}\in H^{1}(\mathbb{A}{},S/P)$.

Any two $c$'s have the same image $H^{1}(\mathbb{Q}{},S^{K}/P^{K})$ because
they have the same image in $H^{1}(\mathbb{A}{},S^{K}/P^{K})$ and we can apply
(\ref{cg6}).

That $a^{-1}(c_{\mathbb{A}})$ is a principal homogeneous space for
$\plim{}C(K)$ follows from Proposition \ref{cg7}.

\subsection{Towards an elementary definition of the fundamental classes}

Are there elementary descriptions of the fundamental classes? By elementary,
we mean involving only class field theory and the cohomology of affine group
schemes. In particular, the definition should not mention abelian varieties,
much less motives.

The definition of $c_{\infty}$ given above is elementary in this sense. But
the definition of $c_{p}$ is not. Wintenberger (1991)\nocite{wintenberger1991}
shows that $\xi_{p}^{K}$ has the following description: choose a prime element
$a$ in $K_{w}$, and let $b=\Nm_{K_{w}/B}(a)$ where $B$ is the maximal
unramified extension of $\mathbb{Q}{}_{p}$ contained in $K_{w}$; define%
\[
\xi^{\prime}(V)=(V\otimes B(\mathbb{F}{}),x\mapsto(1\otimes\sigma
)(bx))\text{.}%
\]
Then $\xi_{p}^{K}\approx\xi^{\prime}$. This gives an elementary description of
$c_{p}^{K}$, but the family $(c_{p}^{K})_{K\in\mathcal{F}{}}$ does not
determine $c_{p}$: there is an exact sequence%
\[
0\rightarrow\plim{}^{1}S^{K}(\mathbb{Q}{}_{p})\rightarrow H^{1}(\mathbb{Q}%
{}_{p},S)\rightarrow\plim{}H^{1}(\mathbb{Q}_{p},{}S^{K})\rightarrow0\text{,}%
\]
but $\varprojlim^{1}S^{K}(\mathbb{Q}{}_{p})\neq0$ because it fails the ML condition.

We shall see later that the Hodge conjecture for CM abelian varieties implies
that there is exactly one distinguished global fundamental class. It seems
doubtful that this can be described elementarily.\pagebreak

\section{Review of gerbes}

In this section, we review some definitions and results from \cite{giraud1971}
and \cite{debremaeker1977a} in the case we need them. Throughout $\mathsf{E}$
is a category with finite fibred products (in particular, a final object $S$)
endowed with a Grothendieck topology. For example, $S$ could be an affine
scheme and $\mathsf{E}$ the category of affine $S$-schemes $\Aff/S$ endowed
with the fpqc topology.

\subsection{Gerbes bound by a sheaf of commutative groups}

\paragraph{Gerbes.}

Recall that a \emph{gerbe on }$S$\emph{\ }is a stack of groupoids
$\varphi\colon\mathsf{G{}}\rightarrow\mathsf{E}$ such that

\begin{enumerate}
\item there exists a covering map $U\rightarrow S$ for which $G{}_{U}$ is nonempty;

\item every two objects of a fibre $\mathsf{G}{}_{U}$ are locally isomorphic
(their inverse images under some covering map $V\rightarrow U$ are isomorphic).
\end{enumerate}

\noindent The gerbe is said to be \emph{neutral }if $\mathsf{G}_{S}$ is nonempty.

Let $\mathbf{x}$ be a cartesian section of $\mathsf{G}/U\rightarrow
\mathsf{E}/U$. Then $\mathcal{A}{}ut(\mathbf{x)}$ is a sheaf of groups on
$\mathsf{E}/U$, which, up to a unique isomorphism, depends only on
$\mathbf{x}(U)$. For $x\in\ob(\mathsf{G}_{U})$, this allows us to define
$\mathcal{A}{}ut(x)=\mathcal{A}{}ut(\mathbf{x})$ with $\mathbf{x}$ any
cartesian section such that $\mathbf{x}(U)=x$.

\paragraph{$A$-gerbes.}

Let $A$ be a sheaf of commutative groups on $\mathsf{E}$. An $A$\emph{-gerbe
on }$S$ is a pair $(\mathsf{G},j)$ where $\mathsf{G}$ is a gerbe on $S$ and
$j(x)$ for $x\in\ob(\mathsf{G})$ is a natural isomorphism $j(x)\colon
\mathcal{A}{}ut(x)\rightarrow A|\varphi(x)$. For example, the gerbe
$\Tors(A)\rightarrow E$ with $\Tors(A)_{U}$ equal to the category of
$A|U$-torsors on $\mathsf{E}/U$ is a (neutral) $A$-gerbe on $S$.

\paragraph{$f$-morphisms.}

Let $f\colon A\rightarrow A^{\prime}$ be a homomorphism of sheaves of
commutative groups, and let $(\mathsf{G},j)$ be an $A$-gerbe on $S$ and
$(\mathsf{G}^{\prime},j^{\prime})$ an $A^{\prime}$-gerbe. An $f$%
\emph{-morphism }from $(\mathsf{G},j)$ to $(\mathsf{G}^{\prime},j^{\prime})$
is a cartesian $E$-functor $\lambda\colon\mathsf{G}\rightarrow\mathsf{G}%
^{\prime}$ such that%
\[
\begin{CD}
\mathcal{A}ut(x) @>{\lambda}>>\mathcal{A}ut(\lambda x)\\
@VVjV@VVj^{\prime}V\\
A|U @>f>>A^{\prime}|U
\end{CD}
\]
commutes for all $U$ and all $x\in\mathsf{G}_{U}$. An $f$-morphism is an
$f$\emph{-equivalence }if it is an equivalence of categories in the usual
sense. If $f$ is an isomorphism, then every $f$-morphism is an $f$%
-equivalence. In particular, every $\id_{A}$-morphism is an $\id_{A}%
$-equivalence (in this context, we shall write $A$-morphism and $A$-equivalence).

Let $\mathsf{G}$ be a trivial $A$-gerbe on $S$. The choice of a cartesian
section $\mathbf{x}$ to $\mathsf{G}\rightarrow\mathsf{E}$ determines an
equivalence of $A$-gerbes $\Hom(\mathbf{x},-)\colon\mathsf{G}\rightarrow
\Tors(A)$. Thus, condition (a) in the definition of a gerbe implies that every
$A$-gerbe is locally isomorphic to $\Tors(A)$.

\paragraph{Morphisms of $f$-morphisms.}

A \emph{morphism }$m\colon\lambda_{1}\rightarrow\lambda_{2}$ of two
$f$-morphisms is simply a morphism of $\mathsf{E}$-functors. With these
definitions, the $A$-gerbes on $S$ form a $2$-category.

Having defined these objects, our next task is to classify them.

\paragraph{Classification of $A$-gerbes.}

One checks easily that $A$-equivalence is an equivalence relation. Giraud
(1971, IV 3.1.1, p247) defines $H^{2}(S,A)$ to be the set of $A$-equivalence
classes of $A$-gerbes, and he then shows that (in the case that $A$ is a sheaf
of \emph{commutative }groups), $H^{2}(S,A)$ is canonically isomorphic to the
usual (derived functor) group (ibid. IV 3.4.2, p261).

\paragraph{Classification of $f$-morphisms and their morphisms.}

Let $f\colon A\rightarrow A^{\prime}$ be a homomorphism of sheaves of
commutative groups, and let $(\mathsf{G},j)$ be an $A$-gerbe and
$(\mathsf{G}^{\prime},j^{\prime})$ an $A^{\prime}$-gerbe on $S$. There is an
$A^{\prime}$-gerbe $\mathsf{HOM}_{f}(\mathsf{G},\mathsf{G}^{\prime})$ on $S$
such that $\mathsf{HOM}_{f}(\mathsf{G},\mathsf{G}^{\prime})_{U}$ is the
category whose objects are the $f$-morphisms $\mathsf{G}|U\rightarrow
\mathsf{G}^{\prime}|U$ and whose morphisms are the morphisms of $f$-morphisms
(\cite{giraud1971}, IV 2.3.2, p218). Its class in $H^{2}(S,A^{\prime})$ is the
difference of the class of $\mathsf{G}^{\prime}$ and the image by $\lambda$ of
the class of $\mathsf{G}$.

\begin{plain}
\label{rg1}We can read off from this the following statements.

\begin{enumerate}
\item There exists an $f$-morphism $\mathsf{G}\rightarrow\mathsf{G}^{\prime}$
if and only if $\lambda$ maps the class of $\mathsf{G}$ in $H^{2}(S,A)$ to the
class of $\mathsf{G}^{\prime}$ in $H^{2}(S,A^{\prime})$ (as $A^{\prime}$ is
commutative, $\mathsf{HOM}_{f}(\mathsf{G},\mathsf{G}^{\prime})$ is neutral if
and only if its class is zero).

\item Let $\lambda_{0}\colon\mathsf{G}\rightarrow\mathsf{G}^{\prime}$ be an
$f$-morphism (assumed to exist). For any other $f$-morphism $\lambda
\colon\mathsf{G}\rightarrow\mathsf{G}^{\prime}$, $\mathcal{H}{}om(\lambda
_{0},\lambda)$ is an $A^{\prime}$-torsor, and the functor $\lambda
\mapsto\mathcal{H}{}om(\lambda_{0},\lambda)$ is an equivalence from the
category whose objects are the $f$-morphisms $\mathsf{G}\rightarrow
\mathsf{G}^{\prime}$ to the category $\Tors(A^{\prime})$. In particular, the
set of isomorphism classes of $f$-morphisms $\mathsf{G}\rightarrow
\mathsf{G}^{\prime}$ is a principal homogeneous space for $H^{1}(S,A^{\prime
})$.

\item Let $\lambda_{1},\lambda_{2}\colon\mathsf{G}\rightarrow\mathsf{G}%
^{\prime}$ be two $f$-morphisms. If they are isomorphic, then the set of
isomorphisms $\lambda_{1}\rightarrow\lambda_{2}$ is a principal homogeneous
space for $H^{0}(S,A^{\prime})=_{\text{df}}A^{\prime}(S)$.
\end{enumerate}
\end{plain}

\begin{exercise}
\label{rg2}Let $S=\Spec k$ with $k$ a field, and let $\mathsf{G}%
\rightarrow\Aff/S$ be a gerbe bound by a separable torus. Let $\bar{S}=\Spec
k^{\text{al}}$, and let $a,b$ be the projection maps $\bar{S}\times_{S}\bar
{S}\rightarrow\bar{S}$. Show that $\mathsf{G}_{\bar{S}}$ is nonempty, and that
for any $x\in\ob\mathsf{G}_{\bar{S}}$, $a^{\ast}x$ and $b^{\ast}x$ are
isomorphic. [Hint: use \ref{cc6}.]
\end{exercise}

\subsection{Gerbes bound by a sheaf of commutative crossed module}

\paragraph{Commutative crossed modules.}

Recall that a \emph{crossed module} is a pair of groups $(A,B)$ together with
an action of $B$ on $A$ and a homomorphism $\rho\colon A\rightarrow B$
respecting this action. It is said to be \emph{commutative} if $A$ and $B$ are
both commutative and the action of $B$ on $A$ is trivial. Thus, a commutative
crossed module is nothing more than a homomorphism of commutative groups.
Similarly, a sheaf of commutative crossed modules is simply a homomorphism
$\rho\colon A\rightarrow B$ of sheaves of commutative groups. A homomorphism
$(f,\phi)\colon(A,B,\rho)\rightarrow(A^{\prime},B^{\prime},\rho^{\prime})$ of
sheaves of commutative crossed modules is a pair of homomorphisms giving a
commutative square%
\[
\begin{CD}
A @>{f}>>A' \\
@VV{\rho}V@VV{\rho^{\prime}}V\\
B@>{\phi}>>B',
\end{CD}
\]
that is, it is a morphism of complexes.

\paragraph{$(A,B)$-gerbes.}

Let $\rho\colon A\rightarrow B$ be a sheaf of commutative crossed modules.
Following \cite{debremaeker1977a}, we define an $(A,B)$-gerbe to be a triple
$(\mathsf{G},\mu,j)$ with $(\mathsf{G},j)$ an $A$-gerb and $\mu$ a $\rho
$-morphism $A\rightarrow\Tors(B)$. For example, $\Tors(A)\overset
{\lambda_{\ast}}{\rightarrow}\Tors(B)$ is an $(A,B)$-gerbe.

\paragraph{$(f,\phi)$-morphisms.}

Let $(f,\phi)\colon(A,B)\rightarrow(A^{\prime},B^{\prime})$ be a homomorphism
of sheaves of commutative crossed modules. Let $(\mathsf{G},\mu,j)$ be an
$(A,B)$-gerbe on $S$ and $(\mathsf{G}^{\prime},\mu^{\prime},j^{\prime})$ an
$(A^{\prime},B^{\prime})$-gerbe. An $(f,\phi)$\emph{-morphism }from
$(\mathsf{G},\mu,j)$ to $(\mathsf{G}^{\prime},\mu^{\prime},j^{\prime})$ is a
pair $(\lambda,i)$ where $\lambda$ is an $f$-morphism $(\mathsf{G}%
,j)\rightarrow(\mathsf{G}^{\prime},j^{\prime})$ and $i$ is an isomorphism of
functors
\[
i\colon\phi_{\ast}\circ\mu\Rightarrow\mu^{\prime}\circ\lambda\colon
\mathsf{G}\rightarrow\Tors(B^{\prime}).
\]
When $(f,\phi)=(\id_{A},\id_{B})$, we speak of an $(A,B)$-morphism
$(\lambda,i)\colon(\mathsf{G},\mu,j)\rightarrow(\mathsf{G}^{\prime}%
,\mu^{\prime},j^{\prime})$.

\paragraph{Morphisms of $(f,\phi)$-morphisms.}

Let $(\lambda_{1},i_{1})$ and $(\lambda_{2},i_{2})$ be $(f,\phi)$-morphims. A
\emph{morphism }$m\colon(\lambda_{1},i_{1})\rightarrow(\lambda_{2},i_{2})$ is
a morphism of $E$-functors $m\colon\lambda_{1}\rightarrow\lambda_{2}$
satisfying the following condition: on applying $\mu^{\prime}$ to $m$, we
obtain morphism of functors $\mu^{\prime}\cdot m\colon\mu^{\prime}\circ
\lambda_{1}\rightarrow\mu^{\prime}\circ\lambda_{2}$; the composite of this
with $i_{1}$ is required to equal $i_{2}$,%
\[
\left(  \mu^{\prime}\cdot m\right)  \circ i_{1}=i_{2}.
\]
With these definitions, the $(A,B)$-gerbes over $S$ form a $2$-category.

Again, we wish to classify these objects.

\paragraph{Classification of $(A,B)$-gerbes.}

An $(A,B)$-morphism is an $(A,B)$-equivalence if $\lambda$ is an equivalence
of categories. Again, $(A,B)$-equivalence is an equivalence relation, and
Debremaeker defines $H^{2}(S,A\rightarrow B)$ to be the set of equivalence
classes. The forgetful functor $(\mathsf{G},\mu,j)\mapsto(\mathsf{G},j)$
defines a map $H^{2}(S,A\rightarrow B)\rightarrow H^{2}(S,A)$.

Because $A\rightarrow B$ is sheaf of \emph{commutative} crossed modules,
$H^{2}(S,A\rightarrow B)$ is in fact canonically isomorphic to the usual
hypercohomology group of the complex $A\rightarrow B$.

\paragraph{Classification of the $(f,\phi)$-morphisms.}

Let $\rho\colon A\rightarrow B$ be a sheaf of commutative crossed modules.
Define an $(A,B)$\emph{-torsor} to be a pair $(P,p)$ with $P$ an $A$-torsor
and $p$ an $S$-point of $\rho_{\ast}P$. A morphism $(P,p)\rightarrow
(P^{\prime},p^{\prime})$ is a morphism $P\rightarrow P^{\prime}$ of
$A$-torsors carrying $p$ to $p^{\prime}$. Define $H^{1}(S,A\rightarrow B)$ to
be the set of isomorphism classes of $(A,B)$-torsors (cf. \cite{deligne1979},
2.4.3; \cite{milne2002c}, \S 1).

Now consider a homomorphism $(f,\phi)\colon(A,B)\rightarrow(A^{\prime
},B^{\prime})$. Let $(\mathsf{G},\mu,j)$ be an $(A,B)$-gerbe over $S$ and
$(\mathsf{G}^{\prime},\mu^{\prime},j^{\prime})$ an $(A^{\prime},B^{\prime}%
)$-gerbe. Assume that there exists an $(f,\phi)$-morphism $(\lambda_{0}%
,i_{0})\colon(\mathsf{G},\mu,j)\rightarrow(\mathsf{G}^{\prime},\mu^{\prime
},j^{\prime})$. As we noted above, the map $\mathcal{\lambda}\mapsto
P(\lambda)=_{\text{df}}\mathcal{H}{}om(\lambda_{0},\lambda)$ is an equivalence
from the category of $f$-morphisms $(\mathsf{G},j)\rightarrow(\mathsf{G}%
^{\prime},j^{\prime})$ to the category of $A^{\prime}$-torsors. From $i$ we
get a point $p(i)\in(\rho_{\ast}^{\prime}P)(S)$, and $(\lambda,i)\mapsto
(P(\lambda),p(i))$ defines an equivalence from the category of $(f,\phi
)$-morphism $(\mathsf{G},\mu,j)\rightarrow(\mathsf{G}^{\prime},\mu^{\prime
},j^{\prime})$ to the category of $(A,B)$-torsors. In particular, we see that
the set of $(f,\phi)$-morphisms $(\mathsf{G},\mu,j)\rightarrow(\mathsf{G}%
^{\prime},\mu^{\prime},j^{\prime})$ is a principal homogeneous space for
$H^{1}(S,A\rightarrow B)$.

\paragraph{Classification of the morphisms of $(f,\phi)$-morphisms.}

Let $(\lambda_{1},i_{1}),(\lambda_{2},i_{2})\colon(\mathsf{G},\mu
,j)\rightarrow(\mathsf{G}^{\prime},\mu^{\prime},j^{\prime})$ be two $(f,\phi
)$-morphisms. If they are isomorphic, then the set of isomorphisms
$(\lambda_{1},i_{1})\rightarrow(\lambda_{2},i_{2})$ is a principal homogeneous
space for
\[
H^{0}(S,A\rightarrow B)=_{\text{df}}\Ker(A(S)\rightarrow B(S)).
\]

\subsection{Gerbes bound by an injective commutative crossed module}

Consider an exact sequence%
\[
0\rightarrow A\overset{\rho}{\rightarrow}B\overset{\sigma}{\rightarrow
}C\rightarrow0
\]
of sheaves of commutative groups.

\paragraph{The group $H^{2}(k,A\rightarrow B)$.}

For any $C$-torsor $P$, there is a gerbe $\mathsf{K}(P)$ whose fibre
$\mathsf{K}(P)_{U}$ over $U$ is the category whose objects are pairs
$(Q,\lambda)$ with $Q$ a $B$-torsor and $\lambda$ a $\rho$-morphism
$Q\rightarrow P$. This has a natural structure of an $A$-gerbe, and the
forgetful functor $(Q,\lambda)\mapsto Q$ endows it with the structure of an
$(A,B)$-gerbe (\cite{debremaeker1977a}).

Conversely, let $(\mathsf{G},\mu,j)$ be an $(A,B)$ gerbe. For any object
$x\in\ob(A_{U})$, $\sigma_{\ast}\mu(x)$ is a $C$-torsor over $U$ endowed with
a canonical descent datum, which gives a $C$-torsor over $S$.

These correspondences define inverse isomorphisms%
\[
H^{2}(k,A\rightarrow B)\cong H^{1}(k,C)\text{.}%
\]

\paragraph{The group $H^{1}(k,A\rightarrow B)$.}

For any point $c\in C(k)$, $\sigma^{-1}(c)$ is an $A$-torsor with $\rho_{\ast
}(\sigma^{-1}(c))=B$. In particular, $\rho_{\ast}(\sigma^{-1}(c))$ has a
canonical point (the identity of $B$), and so $\sigma^{-1}(c)$ has the
structure of an $(A,B)$-torsor.

Conversely, let $(P,p)$ be an $(A,B)$-torsor. For any point $q$ of $P$ in some
covering of $k$, $q-p$ is an element of $B$ whose image in $C$ lies in $C(k)$.

These correspondences define inverse isomorphisms%
\[
H^{1}(k,A\rightarrow B)\cong C(k)\text{.}%
\]

\paragraph{The group $H^{0}(k,A\rightarrow B)$.}

By definition,
\[
H^{0}(k,A\rightarrow B)=\Ker(A(k)\rightarrow B(k))=0.
\]

\newpage

\section{The motivic gerbe}

Grothendieck's construction gives a rigid pseudo-abelian tensor category
$\Mot(\mathbb{F}{})$ of abelian motives over $\mathbb{F}{}$ based on the
abelian varieties over $\mathbb{F}{}$ and using the numerical equivalence
classes of algebraic cycles as correspondences (see, for example,
\cite{saavedra1972}, VI 4). In fact, $\Mot(\mathbb{F}{})$ is Tannakian
(\cite{jannsen1992}). When one assumes the Tate conjecture for abelian
varieties, the fundamental group of $\Mot(\mathbb{F}{})$ becomes identified
with $P$ (e.g., \cite{milne1994a}), and when one assumes Grothendieck's Hodge
standard conjecture, $\Mot(\mathbb{F}{})$ acquires a canonical polarization
(e.g., \cite{saavedra1972}).

Thus, when we assume these two conjectures, $\Mot(\mathbb{F)}$ is a Tannakian
$\mathbb{Q}{}$-category with fundamental group $P$, an $\mathbb{A}_{f}^{p}%
$-valued fibre functor $\omega_{f}^{p}$ (\'{e}tale cohomology), an exact
tensor functor $\omega_{p}\colon\Mot(\mathbb{F}{})\rightarrow\mathsf{R}_{p}$
(crystalline cohomology), and an exact tensor functor $\omega_{\infty}%
\colon\Mot(\mathbb{F}{})\rightarrow\mathsf{R}_{\infty}$ (from the polarization
--- see \cite{deligne1982c}, 5.20). Note that $\omega_{\infty}$ is uniquely
defined up to isomorphism, whereas $\omega^{p}$ and $\omega_{p}$ are uniquely
defined up to a unique isomorphism.

Let $\mathsf{P}=\Mot(\mathbb{F}{})^{\vee}$ be the gerbe of fibre functors on
$\Mot(\mathbb{F}{})$. It is a $P$-gerbe endowed with an object $w^{p}%
=_{\text{df}}\omega_{f}^{p}$ in $\mathsf{P}_{\mathbb{A}_{f}^{p}}$ and
morphisms $w_{p}\colon\mathsf{R}_{p}^{\vee}\rightarrow\mathsf{P}(p)$ and
$w_{\infty}\colon\mathsf{R}_{\infty}^{\vee}\rightarrow\mathsf{P}(\infty)$
(defined by $\omega_{p}$ and $\omega_{\infty}$), where $\mathsf{P}%
(p)=_{\text{df}}\mathsf{P}/\Spec(\mathbb{Q}{}_{p})$ and $\mathsf{P}%
(\infty)=\mathsf{P}/\Spec(\mathbb{R}{})$.

\begin{remark}
Assume the Tate conjecture for all smooth projective varieties over
$\mathbb{F}{}$ and the Hodge standard conjecture for all abelian varieties
over $\mathbb{F}{}$. Then the category of motives based on all smooth
projective varieties is equivalent with its subcategory based on the abelian
varieties (\cite{milne1994a}).\footnote{In particular, it has a canonical
polarization, and the Hodge standard conjecture holds for all smooth
projective varieties in characteristic $p$ (see the arguments in
\cite{milne2002a}).} Therefore, under these assumptions the category of all
motives over $\mathbb{F}{}$ has a system $(\mathsf{P},w^{p},w_{p},w_{\infty})$
attached to it with $\mathsf{P}$ again a $P$-gerbe.
\end{remark}

\subsection{The (pseudo)motivic gerbe}

\begin{plain}
\label{g1}We now drop all assumptions, and consider the existence and
uniqueness of systems $(\mathsf{P},w^{p},w_{p},w_{\infty})$ with

\begin{itemize}
\item $\mathsf{P}$ a $P$-gerbe over $\Spec(\mathbb{Q}),$

\item $w^{p}$ an object of $\mathsf{P}_{\mathbb{A}_{f}^{p}}$,

\item $w_{p}$ an $x_{p}$-morphism $\mathsf{R}_{p}^{\vee}\rightarrow
\mathsf{P}(p),$

\item $w_{\infty}$ an $x_{\infty}$-morphism $\mathsf{R}_{\infty}^{\vee
}\rightarrow\mathsf{P}(\infty).$
\end{itemize}

\noindent Here $x_{p}$ and $x_{\infty}$ are the homomorphisms $\mathbb{G}%
\rightarrow S_{/\mathbb{Q}_{p}}$ and $\mathbb{G}_{m}\rightarrow S_{/\mathbb{R}%
{}}$ described in \S 4.
\end{plain}

\begin{theorem}
\label{g2}There exists a system $(\mathsf{P},w^{p},w_{p},w_{\infty})$ as in
\ref{g1}. If $(\mathsf{P}^{\prime},w^{p\prime},w_{p}^{\prime},w_{\infty
}^{\prime})$ is a second such system, then there exists a $P$-equivalence
$\lambda\colon\mathsf{P}\rightarrow\mathsf{P}^{\prime}$ such that
$\lambda(w^{p})\approx w^{p\prime}$, $\lambda\circ w_{p}\approx w_{p}^{\prime
}$, $\lambda\circ w_{\infty}\approx w_{\infty}^{\prime}$. Any two such
$\lambda$'s are isomorphic by an isomorphism that is unique up to a nonzero
rational number.
\end{theorem}

\begin{proof}
The $P$-equivalence classes of $P$-gerbes $\mathsf{P}$ are classified by
$H^{2}(\mathbb{Q}{},P)$. The condition that there exists a $w^{p}$ is that
$\mathsf{P}/\mathbb{A}_{f}^{p}$ is neutral or, equivalently, that the class of
$\mathsf{P}$ maps to zero in $H^{2}(\mathbb{A}_{f}^{p},P)$; the condition that
there exist a $w_{p}$ is that the class of $\mathsf{P}(p)$ in $H^{2}%
(\mathbb{Q}{}_{p},P)$ is the image of the class of $\mathsf{R}_{p}^{\vee}$
under $x_{p}\colon H^{2}(\mathbb{Q}{}_{p},\mathbb{G}{})\rightarrow
H^{2}(\mathbb{Q}{}_{p},P)$; the condition that there exist a $w_{\infty}$ is
that the class of $\mathsf{P}(\infty)$ in $H^{2}(\mathbb{R}{},P)$ is the image
of the class of $\mathsf{R}_{\infty}^{\vee}$ under $x_{\infty}\colon
H^{2}(\mathbb{R},P)\rightarrow H^{2}(\mathbb{R}{},P)$. Proposition \ref{fc07}
shows that there is a unique class in $H^{2}(\mathbb{Q}{},P)$ with these
properties. Thus, there exists a $P$-gerbe $\mathsf{P}$ for which there exist
$w^{p}$, $w_{p}$, $w_{\infty}$ as in \ref{g1} and if $\mathsf{P}^{\prime}$ is
a second such gerbe, then there exists a $P$-morphism $P\rightarrow P^{\prime
}$.

Now consider two systems $(\mathsf{P},w^{p},w_{p},w_{\infty})$ and
$(\mathsf{P}^{\prime},w^{p\prime},w_{p}^{\prime},w_{\infty}^{\prime})$. The
isomorphism classes of $P$-morphisms $P\rightarrow P^{\prime}$ form a
principal homogeneous space for $H^{1}(\mathbb{Q}{},P)$. On the other hand,
the isomorphism classes of triples $(w^{p},w_{p},w_{\infty})$ for $\mathsf{P}$
(or $\ \mathsf{P}^{\prime}$) form a principle homogeneous space for
$H^{1}(\mathbb{A},P)$. Since the map $H^{1}(\mathbb{Q}{},P)\rightarrow
H^{1}(\mathbb{A}{},P)$ is an isomorphism (\ref{cg5}, \ref{cg5p}), there exists
a $P$-morphism $\lambda\colon\mathsf{P}\rightarrow\mathsf{P}^{\prime}$
carrying the isomorphism class of $(w^{p},w_{p},w_{\infty})$ into that of
$(w^{p\prime},w_{p}^{\prime},w_{\infty}^{\prime})$, and it is unique up to
isomorphism. If $\lambda_{1}$ and $\lambda_{2}$ are isomorphic $P$-morphisms
$P\rightarrow P^{\prime}$, then the set of isomorphisms $\lambda
_{1}\rightarrow\lambda_{2}$ is a principal homogeneous space for
$P(\mathbb{Q}{})$, which equals $\mathbb{Q}{}^{\times}$ (\ref{cg5r}).
\end{proof}

\begin{remark}
\label{gc2}When we replace $\Mot(\mathbb{F}{})$ with its subcategory
$\Mot_{0}(\mathbb{F}{})$ of motives of weight $0$, we obtain similar results,
except that its gerbe $\mathsf{P}_{0}$ is a $P_{0}$-gerbe. Theorem \ref{g2}
holds with the only change being that now any two $\lambda$'s are uniquely
isomorphic (because $P_{0}(\mathbb{Q}{})=0$).
\end{remark}

\subsection{The pseudomotivic groupoid}

Let $\nu_{p}$ be the forgetful $B(\mathbb{F}{})$-valued fibre functor on the
Tannakian category $\mathsf{R}_{p}$, and let $\mathfrak{G}{}_{p}%
=\underline{\Aut}^{\otimes}(\nu_{p})$. It is a $\mathbb{Q}{}_{p}^{\text{al}%
}/\mathbb{Q}{}_{p}$-groupoid with kernel $\mathfrak{G}{}_{p}^{\Delta
}=\mathbb{G}{}$. The fibre functor $v_{p}$ is an object of $(\mathsf{R}%
_{p}^{\vee})_{B(\mathbb{F}{})}$, and $\mathfrak{G}{}_{p}$ can also be
described as the groupoid of automorphisms of this object (in the sense of
\cite{deligne1990}, 3.4).

Let $\nu_{\infty}$ be the forgetful $\mathbb{C}$-valued fibre functor on the
Tannakian category $\mathsf{R}_{\infty}$, and let $\mathfrak{G}{}_{\infty
}=\underline{\Aut}^{\otimes}(\nu_{\infty})$. It is an $\mathbb{C}/\mathbb{R}%
$-groupoid with kernel $\mathfrak{G}{}_{\infty}^{\Delta}=\mathbb{G}%
_{m/\mathbb{C}{}}$. The fibre functor $\nu_{\infty}$ is an object of
$(\mathsf{R}_{\infty}^{\vee})_{\mathbb{C}{}},$ and $\mathfrak{G}{}_{\infty}$
can also be described as the groupoid of automorphisms of this object.

For $l\neq p,\infty$, let $\mathfrak{G}{}_{l}$ be the trivial $\mathbb{Q}%
{}_{l}^{\text{al}}/\mathbb{Q}{}_{l}$-groupoid.

Let $(\mathsf{P},w^{p},w_{p},w_{\infty})$ be as in (\ref{g1}). Then, because
$H^{i}(\mathbb{Q}^{\text{al}},P)=0$ for $i>0$ (\ref{cc6}), there exists an
$x\in\ob(\mathsf{P}_{\mathbb{Q}{}^{\text{al}}})$, and any two such objects are
isomorphic. Let $S=\Spec\mathbb{Q}{}$ and $\bar{S}=\Spec\mathbb{Q}%
{}^{\text{al}}$. Let $\mathfrak{P}$ be the $\bar{S}/S$-groupoid of
automorphisms of $x$: for any $\bar{S}\times_{S}\bar{S}$-scheme $(b,a)\colon
T\rightarrow\bar{S}\times_{S}\bar{S}$, $\mathfrak{P}{}(T)$ is the set of
isomorphisms $a^{\ast}x\rightarrow b^{\ast}x$. Its kernel $\mathfrak{P}%
{}^{\Delta}=P$. For $l\neq p,\infty$, $w^{p}$ defines a homomorphism
$\zeta_{l}\colon\mathfrak{G}{}_{l}\rightarrow\mathfrak{P}(l)$ where
$\mathfrak{P}{}(l)$ is the $\mathbb{Q}{}_{l}^{\text{al}}/\mathbb{Q}{}_{l}%
$-groupoid obtained from $\mathfrak{P}$ by base change. Moreover, $w_{p}$
defines a homomorphism $\zeta_{p}\colon\mathfrak{G}_{p}\rightarrow
\mathfrak{P}{}(p)$ and $w_{\infty}$ defines a homomorphism $\zeta_{\infty
}\colon\mathfrak{G}_{\infty}\rightarrow\mathfrak{P}{}(\infty)$.

\begin{proposition}
\label{gc4}

\begin{enumerate}
\item The system $(\mathfrak{\mathfrak{P}{}},(\zeta_{l})_{l})$ satisfies the
following conditions:

\begin{enumerate}
\item $(\mathfrak{P}{}^{\Delta},\zeta_{p}^{\Delta},\zeta_{\infty}^{\Delta
})=(P,x_{p},x_{\infty});$

\item the morphisms $\zeta_{l}$ for $l\neq p,\infty$ are induced by a section
of $\mathfrak{P}{}$ over $\Spec(\overline{\mathbb{A}_{f}^{p}}\otimes
_{\mathbb{A}_{f}^{p}}\overline{\mathbb{A}_{f}^{p}})$ where $\overline
{\mathbb{A}_{f}^{p}}$ is the image of the map $\mathbb{\bar{Q}}\otimes
_{\mathbb{Q}{}}\mathbb{A}_{f}^{p}\rightarrow%
{\textstyle\prod_{l\neq p,\infty}}
\mathbb{Q}{}_{\ell}^{\text{al}}{}.$
\end{enumerate}

\item Let $(\mathfrak{P}^{\prime},(\zeta_{l}^{\prime}))$ be the system
attached to $x^{\prime}\in\ob(\mathsf{P}_{\mathbb{Q}{}^{\text{al}}}^{\prime}%
)$. Then there exists an isomorphism $\alpha\colon\mathfrak{P}\rightarrow
\mathfrak{P}^{\prime}$ such that, for all $l$, $\zeta_{l}^{\prime}$ is
isomorphic to $\alpha\circ\zeta_{l}$, and any two such $\alpha$'s are isomorphic.
\end{enumerate}
\end{proposition}

\begin{proof}
Straightforward consequence of the theorem.
\end{proof}

\begin{definition}
Any system $(\mathfrak{\mathfrak{P}{},(\zeta}_{l})_{l})$ arising from a system
$(\mathsf{P},w^{p},w_{p},w_{\infty})$ as in (\ref{g1}) and an object $x\in
P_{\mathbb{Q}{}^{\text{al}}}$ will be called a \emph{pseudomotivic groupoid}.
\end{definition}

\paragraph{Notes.}

In \cite{milne1992a}, 3.27, a pseudomotivic groupoid is defined to be any
system $(\mathfrak{P}{},(\zeta_{l}))$ satisfying condition (\ref{gc4}a(i)).
Theorem 3.28 (ibid.) then states: there exists a pseudomotivic groupoid
$(\mathfrak{P}{},(\zeta_{l}))$; if $(\mathfrak{P}{},(\zeta_{l}))$ is a second
pseudomotivic groupoid, then there is an isomorphism $\alpha\colon
\mathfrak{P}{}\rightarrow\mathfrak{P}{}^{\prime}$ such that $\zeta_{l}%
^{\prime}\approx\alpha\circ\zeta_{l}$, and $\alpha$ is uniquely determined up
to isomorphism. Only a brief indication of proof is given, and the theorem is
credited to Langlands and Rapoport (1987). A sketch of a proof is given in
\cite{milne1994a}, 3.31.

As Reimann points out (1997, p120), I inadvertently omitted the condition
(\ref{gc4}a(ii)) from the definition of the pseudomotivic
groupoid.\footnote{Condition (\ref{gc4}a(ii)) is needed for the statement of
the conjecture of Langlands and Rapoport.} Moreover, the argument sketched in
\cite{milne1994a} (based on that in \cite{langlands1987}) shows only that
there exists an isomorphism $\alpha$ for which $\zeta_{l}^{\prime}$ is
algebraically isomorphic to $\alpha\circ\zeta_{l}$ (i.e., becomes isomorphic
when projected to an algebraic quotient).

However, there is a more serious criticism of the four articles just cited
(and others), namely, in each article $H^{i}(k,P)$ is taken to be $\varprojlim
H^{i}(k,P^{K})$ instead of the fpqc group $H^{i}(k,P)$ which is, in fact, the
group that classifies the various objects.\footnote{In \cite{milne1994a},
footnote p441, I correctly note that the gerbes are classified by the fpqc
group $H^{2}(k,P)$, but then claim that this group equals $\varprojlim
H^{2}(k,P^{K})$, giving as a reference \cite{saavedra1972}, III 3.1. In fact,
Saavedra proves results only for algebraic groups. I don't know why I thought
the results held for affine group schemes, except perhaps I confused the
(true) statement that cohomology commutes with products with the (false)
statement that it commutes with inverse limits.} This amounts to ignoring the
terms $\varprojlim^{1}H^{i-1}(k,P^{K})$, which are not all zero. Thus, some of
the proofs in these papers are inadequate.

The papers cited above all work with groupoids.\footnote{This is obscured in
\cite{langlands1987} by the authors' calling their groupoids gerbes.} Here, I
have preferred to work with gerbes because their attachment to Tannakian
categories is canonical (the groupoid of a Tannakian category depends on the
choice of a fibre functor), are more directly related to nonabelian
cohomology, and are, in some respects, easier.

\newpage

\section{The motivic morphism of gerbes}

Assume the Hodge conjecture for complex abelian varieties of CM-type. This
implies the Tate conjecture and the Hodge standard conjecture for abelian
varieties over finite fields (\cite{milne1999b}; \cite{milne2001a}). Hence we
get a category of abelian motives $\Mot(\mathbb{F}{})$ with fibre functors as
in the last section. Moreover, the category $\CM(\mathbb{Q}^{\text{al}})$ of
CM-motives over $\mathbb{Q}{}^{\text{al}}$, defined using Deligne's Hodge
classes, coincides with that defined using algebraic cycles modulo numerical
equivalence. Thus, we get a reduction functor $R\colon\CM(\mathbb{Q}%
^{\text{al}})\rightarrow\Mot(\mathbb{F}{})$ (exact tensor functor of
$\mathbb{Q}{}$-categories), and hence a morphism of gerbes of fibre functors%
\[
R^{\vee}\colon\Mot(\mathbb{F}{})^{\vee}\rightarrow\CM(\mathbb{Q}^{\text{al}%
})^{\vee}\text{.}%
\]
It follows from Shimura-Taniyama theory that this is a $\rho$-morphism, where
$\rho\colon P\rightarrow S$ is the homomorphism defined in \S 3.

The Betti fibre functor $\omega_{B}$ on $\CM(\mathbb{Q}^{\text{al}})$ defines
an $S$-equivalence of gerbes
\[
\omega\mapsto\mathcal{H}{}om(\omega_{B},\omega)\colon\CM(\mathbb{Q}%
^{\text{al}})^{\vee}\rightarrow\Tors(S)\text{.}%
\]
On composing $R^{\vee}$ with this, we obtain a $\rho$-morphism $\mu
\colon\Mot(\mathbb{F}{})^{\vee}\rightarrow\Tors(S)$, i.e., a $(P,S)$-gerbe.
For all $l\neq p,\infty$, $\mu(\omega_{l})$ is a trivial $S_{\mathbb{Q}{}_{l}%
}$-torsor${}$. Moreover, the composites of $w_{p}$ and $w_{\infty}$ with $\mu$
are isomorphic to the functors $\xi_{p}^{\vee}$ and $\xi_{\infty}^{\vee}$,
where $\xi_{p}$ and $\xi_{\infty}$ are as in \S 4.

\subsection{The (pseudo)motivic morphism of gerbes}

\begin{plain}
\label{g3}We now drop all assumptions, and consider the existence and
uniqueness of systems $(\mathsf{P},\mu,w^{p},w_{p},w_{\infty})$ with

\begin{itemize}
\item $(\mathsf{P},\mu)$ a $(P,S)$-gerbe,

\item $w^{p}$ an object of $\mathsf{P}_{\mathbb{A}_{f}^{p}}$ such that
$\mu(w^{p})$ is a trivial $S_{\mathbb{A}_{f}^{p}}$-torsor;

\item $w_{p}$ an $x_{p}$-morphism $\mathsf{R}_{p}^{\vee}\rightarrow
\mathsf{P}(p)$ such that $\mu\circ w_{p}\approx\xi_{p}^{\vee}$;

\item $w_{\infty}$ an $x_{\infty}$-morphism $\mathsf{R}_{\infty}^{\vee
}\rightarrow\mathsf{P}(\infty)$ such that $\mu\circ w_{\infty}\approx
\xi_{\infty}^{\vee}$.
\end{itemize}
\end{plain}

\begin{theorem}
\label{mg1}The set of $(P,S)$-equivalence classes of $(P,S)$-gerbes
$(\mathsf{P},\mu)$ that can be completed to a system $(\mathsf{P},\mu
,w^{p},w_{p},w_{\infty})$ as in (\ref{g3}) is a $\varprojlim_{\mathcal{F}{}%
}C(K)$-principal homogeneous space, where $\mathcal{F}{}$ is the set of
CM-subfields of $\mathbb{Q}{}^{\text{al}}$ of finite degree over $\mathbb{Q}%
{}$; in particular, it is nonempty.
\end{theorem}

\begin{proof}
The $(P,S)$-gerbes $(\mathsf{P},\mu)$ are classified by $H^{2}(\mathbb{Q}%
{},P\rightarrow S)\cong H^{1}(\mathbb{Q}{},S/P)$. The condition that there
exists a $w^{p}$ with $\mu(w^{p})$ neutral is that $(\mathsf{P},\mu)$ becomes
equivalent with the trivial $(P,S)$-gerbe $\Tors(P)\overset{\rho_{\ast}%
}{\rightarrow}\Tors(S)$ over $\mathbb{A}_{f}^{p}$, i.e., that the class of
$\ (\mathsf{P},\mu)$ maps to zero in $H^{1}(\mathbb{A}_{f}^{p},S/P)$; the
condition that there exists a $w_{p}$ with $\mu\circ w_{p}\approx\xi_{p}%
^{\vee}$ is that the class of $(\mathsf{P},\mu)$ over $\mathbb{Q}{}_{p}$ is
the local fundamental class $c_{p}$; the condition that there exists a
$w_{\infty}$ with $\mu\circ w_{\infty}\approx\xi_{\infty}^{\vee}$ is that the
class of $(\mathsf{P},\mu)$ over $\mathbb{R}$ is the local fundamental class
$c_{\infty}$. Theorem \ref{fc02} shows that there exists such a class in
$H^{1}(\mathbb{Q}{},S/P)$, and that the set of them is a principal homogeneous
space under $\varprojlim C(K)$. This completes the proof.
\end{proof}

\begin{remark}
\label{g5}If $(\mathsf{P},\mu)$ and $(\mathsf{P}^{\prime},\mu^{\prime})$ are
$(P,S)$-gerbes that can be completed to systems as in (\ref{g3}), then, for
all $K\in\mathcal{F}{}$, $(\mathsf{P}^{K},\mu^{K})$ and $(\mathsf{P}^{\prime
K},\mu^{\prime K})$ are $(P^{K},S^{K})$-equivalent $(P^{K},S^{K})$-gerbes. The
point is that any two fundamental classes have the same image in
$H^{1}(\mathbb{Q}{},S^{K}/P^{K})$ (\ref{fc02}).
\end{remark}

\begin{remark}
\label{g4}Let $(\mathsf{P}^{\prime},\mu^{\prime},w^{p\prime},w_{p}^{\prime
},w_{\infty}^{\prime})$ be a second system as in \ref{g3}. Then, even if
$(\mathsf{P},\mu)$ and $(\mathsf{P}^{\prime},\mu^{\prime})$ are $(P,S)$%
-equivalent, there may not be a $(P,S)$-equivalence $(\lambda,i)\colon
(\mathsf{P},\mu)\rightarrow(\mathsf{P}^{\prime},\mu^{\prime})$ such that
$\lambda(w^{p})\approx w^{p\prime}$, $\lambda\circ w_{p}\approx w_{p}^{\prime
}$, $\lambda\circ w_{\infty}\approx w_{\infty}^{\prime}$. The point is that
the set of isomorphism classes of $(P,S)$-equivalences $(\mathsf{P}%
,\mu)\rightarrow(\mathsf{P}^{\prime},\mu^{\prime})$ (if nonempty) is a
principle homogeneous space under $(S/P)(\mathbb{Q}{})$, while the set of
isomorphism classes of triples $(w^{p},w_{p},w_{\infty})$ satisfying
(\ref{g3}) for $(P,\mu)$ is a principal homogeneous space under
$(S/P)(\mathbb{A}{})$, and $(S/P)(\mathbb{Q}{})\rightarrow(S/P)(\mathbb{A}{})$
is not surjective (cf. \ref{cg2}).
\end{remark}

\begin{remark}
\label{g6}Let $(\mathsf{P},\mu,w^{p},w_{p},w_{\infty})$ be as in \ref{g3}. On
choosing an object $x\in\mathsf{P}_{\mathbb{Q}^{\text{al}}}$ and an
isomorphism of $\mu(x)$ with the trivial torsor, we obtain a morphism of
groupoids $\mathfrak{P}{}\rightarrow\mathfrak{G}{}_{S}$ with $(\mathfrak{P}%
{}\rightarrow\mathfrak{G}{}_{S})^{\Delta}\cong(P\overset{\rho}{\rightarrow}%
S)$. Following \cite{pfau1993}, we define the quasimotivic groupoid to be
$\mathfrak{P}{}\times_{\mathfrak{G}{}_{S}}\mathfrak{G}{}_{T}$ where
$T=\varprojlim(\mathbb{G}_{m})_{L/\mathbb{Q}{}}$ (limit over all subfields $L$
of $\mathbb{Q}{}^{\text{al}}$ of finite degree over $\mathbb{Q}{}$).
\end{remark}

\begin{remark}
\label{mg7}To state the conjecture of Langlands and Rapoport
(\cite{langlands1987}, p169) for Shimura varieties $\Sh(G,X)$ with rational
weight and $G^{\text{der}}$ simply connected, only the pseudomotivic groupoid
is needed. For Shimura varieties with rational weight and $G^{\text{der}}$ not
necessarily simply connected, one needs\footnote{Reimann (1997) overlooks
this: in his Conjecture B3.7 it is necessary to require that $G^{\text{der}}$
be simply connnected.} the morphism $\mathfrak{P}{}\rightarrow\mathfrak{G}%
{}_{S}$ (\cite{milne1992a}). For arbitrary Shimura varieties, one needs a
quasimotivic groupoid.
\end{remark}

\begin{remark}
\label{mg6}Assume

\begin{enumerate}
\item motivated classes (in the sense of \cite{andre1996}) reduce modulo $w$
to motivated classes;

\item the intersection number of any two motivated classes of complementary
dimension on a smooth projective variety over $\mathbb{F}{}$ is a rational
number.\footnote{Only the second of these statements appears to be beyond
reach at present.}
\end{enumerate}

\noindent Then the arguments of \cite{milne1999b} and \cite{milne2002a} show
that the Tate conjecture and the Hodge standard conjecture hold for abelian
varieties over $\mathbb{F}{}$ with \textquotedblleft algebraic
class\textquotedblright\ replaced with \textquotedblleft motivated
class\textquotedblright.\footnote{The Hodge standard conjecture then holds for
algebraic classes on abelian varieties in prime characteristic!} Thus, under
the assumption of (a) and (b), there exists a Tannakian category of abelian
motives $\Mot(\mathbb{F}{})$ over $\mathbb{F}{}$ (defined using motivated
correspondences) and a reduction functor $\CM(\mathbb{Q}^{\text{al}%
})\rightarrow\Mot(\mathbb{F}{})$. From this, we obtain a \emph{well-defined}
system $(P,\mu,w^{p},w_{p},w_{\infty})$ as in \ref{g3}.
\end{remark}

\newpage

\bibliographystyle{amsxport}
\bibliography{refs}

\end{document}